\documentclass{amsart}
\usepackage[utf8]{inputenc}

\usepackage{float}
\usepackage{amsmath}%
\usepackage{amsfonts}%
\usepackage{amsthm}
\usepackage{amssymb}%
\usepackage{graphicx}
\usepackage{hyperref}
\usepackage{csquotes}
\usepackage[all,cmtip]{xy}
\usepackage{geometry}
\theoremstyle{plain}
\newtheorem{theorem}{Theorem}[section]
\newtheorem{lemma}[theorem]{Lemma}

\newtheorem{proposition}[theorem]{Proposition}

\makeatletter

\newenvironment{proofofprop}[1]{\par
  \pushQED{\qed}%
  \normalfont \topsep6\p@\@plus6\p@\relax
  \trivlist
  \item[\hskip\labelsep
        \emph{Proof of Proposition \ref{vanishingprop}.}]\ignorespaces
}{%
  \popQED\endtrivlist\@endpefalse
}
\makeatother
\theoremstyle{definition}
\newtheorem{choice}[theorem]{Choice}

\newtheorem{definition}[theorem]{Definition}

\theoremstyle{remark}
\newtheorem{remark}[theorem]{Remark}

\newcommand{\sqcommdiag}[8]{\[
\xymatrixcolsep{4pc}\xymatrix{
#1 \ar[r]^{#6} \ar[d]_{#5} & #2 \ar[d]^{#7} \\
 #3 \ar[r]_{#8}& #4
}
\]
}

\newcommand{\Z}{\mathbb{Z}}
\newcommand{\Q}{\mathbb{Q}}
\newcommand{\R}{\mathbb{R}}

\newcommand{\spec}{\mathrm{Spec}}

\renewcommand{\P}{\mathbb{P}}
\newcommand{\F}{\mathbb{F}}

\newcommand{\al}{\alpha}

\newcommand{\dd}{\mathrm{d}}

\begin{document}
\title{Quadratic points on modular curves with infinite Mordell--Weil group}
\author{Josha Box}
\begin{abstract}
Bruin and Najman \cite{bruin} and Ozman and Siksek \cite{ozman} have recently determined the quadratic points on each modular curve $X_0(N)$ of genus 2, 3, 4, or 5 whose Mordell--Weil group has rank 0. In this paper we do the same for the $X_0(N)$ of genus 2, 3, 4, and 5 and positive Mordell--Weil rank. The values of $N$ are 37, 43, 53, 61, 57, 65, 67 and 73.

The main tool used is a relative symmetric Chabauty method, in combination with the Mordell--Weil sieve. Often the quadratic points are not finite, as the degree 2 map $X_0(N)\to X_0(N)^+$ can be a source of infinitely many such points. In such cases, we describe this map and the rational points on $X_0(N)^+$, and we specify the exceptional quadratic points on $X_0(N)$ not coming from $X_0(N)^+$. In particular we determine the $j$-invariants of the corresponding elliptic curves and whether they are $\Q$-curves or have complex multiplication.
\end{abstract}
\date{\today}
\thanks{
\indent During the work on this article, the author was supported by an EPSRC DTP studentship.}
\maketitle
\tableofcontents

\section{Introduction}
%
Let $N$ be a positive integer. In a celebrated paper \cite{mazur1}, Mazur determined exactly which modular curves $X_1(N)$ admit non-cuspidal rational points and, for prime values of $N$, he repeated this for $X_0(N)$ \cite{mazur0}. These results for $X_0(N)$ were later extended to composite levels by Kenku \cite{kenku}. Since this work of Mazur, people have been interested in obtaining similar results for points of low degree $d$. 

For $d=2$ (Kamienny \cite{kamienny}), $d=3$ (Parent \cite{parent}), and $d\in \{4,5,6\}$ (Derickx, Kamienny, Stein and Stoll \cite{dkss}), the prime values of $N$ such that $X_1(N)$ has a non-cuspidal degree $d$ point have been determined explicitly. Moreover, Merel's uniform boundedness theorem \cite{merel} proves the existence of an upper bound $B_d$ such that $X_1(N)$ has no non-cuspidal points for primes $N\geq B_d$. The low degree points on $X_0(N)$, however, are naturally more abundant than on $X_1(N)$, which complicates their study. 

There are two obvious potential sources of infinitely many quadratic points on $X_0(N)$: a degree 2 map $X_0(N)\to \P^1$ over $\Q$ (which exists if and only if $X_0(N)$ is hyperelliptic) and a degree 2 map over $\Q$ to an elliptic curve with infinite Mordell--Weil group (the existence of which implies that $X_0(N)$ is bielliptic). In both of these cases the  rational points on the image give rise to infinitely many quadratic points on $X_0(N)$. Using Faltings' theorem \cite{faltings} on abelian varieties, Abramovich and Harris \cite{abramovich} show that the set of quadratic points on $X_0(N)$ is infinite in these two cases \emph{only}. Building on the work of Harris and Silverman \cite{harris}, Bars \cite{bars} then determined all bielliptic $X_0(N)$ and decided that exactly 10 of these have a quotient of positive Mordell--Weil rank. Moreover, Ogg \cite{ogg} decided for which 19 values of $N$ the curve $X_0(N)$ is hyperelliptic. Consequently, the values of $N$ where $X_0(N)$ has genus at least 2 and admits infinitely many quadratic points are  22, 23, 26, 28, 29, 30, 31, 33, 35, 37, 39, 40, 41, 43, 46, 47, 48, 50, 53, 59, 61, 65, 71, 79, 83, 89, 101 and 131. (Of course there are also infinitely many quadratic points on each of the genus 0 and genus 1 curves.) The careful reader will have noticed there must be one curve that is both hyperelliptic and bielliptic with quotient of positive rank: this is the infamous $X_0(37)$.

 Even when infinite, the quadratic points on $X_0(N)$ can be described. Recently, Bruin and Najman \cite{bruin} determined the finitely many exceptional quadratic points (those not coming from $\P^1(\Q)$) on all hyperelliptic $X_0(N)$ except for $N=37$, and they moreover proved in those cases that the quadratic points coming from $\P^1(\Q)$ correspond to $\Q$-curves. The remaining case $N=37$ is also the only hyperelliptic one where $J_0(N)(\Q)$ is infinite. Subsequently, Ozman and Siksek \cite{ozman} determined all finitely many quadratic points on $X_0(N)$ when it is non-hyperelliptic of genus 2, 3, 4 or 5 and the Mordell--Weil group $J_0(N)(\Q)$ is finite.\\

In this paper, we complement the work of Ozman and Siksek by describing all quadratic points on the modular curves $X_0(N)$ of genus 2, 3, 4 and 5 whose Mordell--Weil group is infinite. The corresponding values of $N$, together with the genus $g$ of $X_0(N)$ and the rank $r$ of $J_0(N)(\Q)$ are listed in the table below.
\vspace{-0.35cm}%
\begin{table}[h!]
    \centering 
    \begin{tabular}{c c c c c c c c c }
        $N$ & 37 & 43 & 53 & 61 & 57 & 65 & 67 & 73   \\  \hline
        $g$ & 2 & 3 & 4 & 4 & 5 & 5 & 5 & 5\\[-1.2ex]
        $r$ &1 & 1 & 1 & 1 & 1 & 1 & 2& 2\\ 
       \end{tabular}
    \label{Nlist}
\end{table}%
\vspace{-0.35cm}%
The values of $N$ such that $X_0(N)$ has genus 2, 3, 4 or 5 can be determined via the genus formula, and the rank can be found using the Modular Arithmetic Geometry package in \texttt{Magma} that uses Stein's algorithms \cite{stein}, \cite{stein2}; this was done in Section 2 of \cite{ozman}. 

We denote by $X_0(N)^+$ the quotient of $X_0(N)$ by the Atkin--Lehner involution $w_N$ corresponding to $N$. As mentioned before, the case $N=37$ is special because $X_0(N)$ is hyperelliptic. The Atkin--Lehner involution and the hyperelliptic involution do not coincide, causing both $\P^1$ and the rank 1 elliptic curve $X_0(37)^+$ to contribute infinitely many quadratic points. Despite this, a description of all quadratic points is still possible, albeit slightly less satisfying than in the other cases. We give this description in Proposition \ref{prop37} in the final Section \ref{hypersection}.

Until then, we assume $N$ to be in the above list but unequal to 37.  The main trick to studying quadratic points on a curve $X$ is to instead study the rational points on its symmetric square $X^{(2)}$. This set $X^{(2)}(\Q)$ consist of pairs $\{P,\overline{P}\}$ of a genuinely quadratic point and its conjugate, as well as pairs $\{P,Q\}$ of rational points $P,Q\in X(\Q)$. Following the suggestion of Ozman and Siksek in \cite{ozman}, we have studied the rational points on $X_0(N)^{(2)}$ using the (relative) symmetric Chabauty method developed by Siksek in \cite{siksek}, in combination with a Mordell--Weil sieve. We describe the Chabauty method in more detail in Sections \ref{coleman}--\ref{relativechab}, while the sieve is explained in Section \ref{sieve}.

The values of $N$ in our list for which $X_0(N)$ has finitely many quadratic points are 57, 67 and 73. The symmetric Chabauty method  was successful in determining all twelve quadratic points on $X_0(57)$. For the remaining values of $N$ we used the relative symmetric Chabauty method also developed by Siksek in the same paper \cite{siksek}. In those cases we found that $\mathrm{rk}\,J_0^+(N)(\Q)=\mathrm{rk}\,J_0(N)(\Q)$, allowing us to use Chabauty relative to the degree 2 map $X_0(N)\to X_0(N)^+$ to determine all finitely many quadratic points on $X_0(N)$ not coming from $X_0(N)^+(\Q)$. For $N\in \{67,73\}$, it then remains to determine the rational points on $X_0(N)^+$, which is hyperelliptic of genus 2 and rank 2 in both cases. These curves $X_0(67)^+$ and $X_0(73)^+$ turn out to be beautiful test cases for the explicit quadratic Chabauty method developed by Balakrishnan and Dogra \cite{dogra} and Balakrishnan, Dogra, M\"uller, Tuitman and Vonk \cite{tuitman} following Kim's work \cite{kim1}, \cite{kim2} on non-abelian Chabauty. Their rational points were determined recently in joint work by Balakrishnan, Best, Bianchi, Lawrence, M\"uller,  Triantafillou and Vonk \cite{bbbmtv}, thus giving us a complete list of all quadratic points on $X_0(67)$ and $X_0(73)$. For $N\in \{43,53,61,65\}$, $X_0(N)^+$ is a rank 1 elliptic curve, and we have computed explicitly the map $X_0(N)\to X_0(N)^+$ as well as generators for $X_0(N)^+(\Q)$. The \texttt{Magma} code to verify all these computations can be found at
\[
\texttt{\href{https://github.com/joshabox/quadraticpoints/}{https://github.com/joshabox/quadraticpoints/}} \;.
\]

\begin{theorem}[Main theorem]
All finitely many quadratic points on the modular curves $X_0(N)$ for $N\in \{57,67,73\}$ are as described in the tables in Section \ref{results}. For $N\in \{43,53,61,65\}$, the set of quadratic points on $X_0(N)$ is infinite. The tables in Section \ref{results} give all the finitely many quadratic points not coming from $X_0(N)^+(\Q)$ as well as generators for the Mordell--Weil groups $X_0(N)^+(\Q)$ of the elliptic curves $X_0(N)^+$.
\end{theorem}
We call points on $X_0(N)$ not coming from $X_0(N)^+(\Q)$ \emph{exceptional}. 

If a quadratic point is non-exceptional, the Atkin--Lehner involution $w_N$ defines an $N$-isogeny from the corresponding elliptic curve to its conjugate. This implies that this elliptic curve is a $\Q$-curve: it is $\overline{\Q}$-isogenous to all its Galois conjugates. In attempts to prove modularity results over quadratic fields one is often  led to the study of quadratic points on modular curves; see for example the proof of Freitas, Le Hung and Siksek \cite{freitas} for elliptic curves over totally real quadratic fields. Since $\Q$-curves are automatically modular \cite{ribet}, it is of interest to know the remaining quadratic points.\\
%
%

The author would like to express his sincere gratitude to Samir Siksek for various invaluable suggestions and enlightening conversations, which have without doubt improved this work greatly. 

We also thank the anonymous referees for their detailed feedback including several useful suggestions and corrections that improved the quality of this paper.

\section{Relative Symmetric Chabauty and the Mordell--Weil sieve}\label{chabauty} In this section we describe the relative symmetric Chabauty method developed by Siksek in \cite{siksek}. Thoughout this Section, let $X/\Q$ be a (smooth projective) non-hyperelliptic curve of genus $g\geq 3$ with Jacobian $J$. Let $p$ be a prime of good reduction for $X$.

\subsection{Chabauty--Coleman}\label{coleman} When the rank $r$ of $J(\Q)$ and the genus $g$ of $X$ satisfy the \emph{Chabauty assumption}
\[
r<g,
\]
the method of Chabauty \cite{chabauty} and Coleman \cite{coleman} can be used to find a finite set of points containing $X(\Q_p)\cap \overline{J(\Q)}$, where the closure is inside the $p$-adic topology on $J(\Q_p)$. Here we use a rational point $\infty$ to embed $X\to J$ via $P\mapsto [P-\infty]$.  In order to bound $X(\Q_p)\cap \overline{J(\Q)}$, one needs to find at least one global differential $\omega \in \Omega_{X/\Q_p}(X)$ whose Coleman integrals $\int_0^D\omega$ vanish on all $D\in X(\Q_p)\cap J(\Q)$. We call such a differential a \emph{vanishing differential}.  Such integrals can locally be written as a $p$-adic power series in a local parameter and one can then bound its number of zeros using the theory of Newton polygons. Combining this information for several primes $p$, one hopes to determine all the rational points on $X$. A great source for more information about Chabauty--Coleman is the survey article by McCallum and Poonen \cite{poonen}.

\subsection{Symmetric Chabauty} \label{symmchab}
\subsubsection{Introduction}To study all quadratic points on $X$ at once, one often  instead studies the rational points on the symmetric square $X^{(2)}$. We denote points on $X^{(2)}$ as a 2-set $\{P,Q\}$ of points $P,Q\in X$. Recall that the set $X^{(2)}(\Q)$ consists of pairs $\{P,\overline{P}\}$ of a genuinely quadratic point $P$ on $X$ and its Galois conjugate, as well as pairs $\{P,Q\}$ of rational points $P,Q$ on $X$. As $X$ is non-hyperelliptic of genus $g\geq 3$, also $X^{(2)}$ can be embedded in $J$ via $\{P,Q\}\mapsto [P+Q-2\cdot \infty]$. Now one can use the same method to potentially determine a set containing
\[
X^{(2)}(\Q_p)\cap \overline{J(\Q)},
\]
where $p$ is a prime. However, as $X^{(2)}$ is 2-dimensional, we now need at least two linearly independent differentials vanishing on $X^{(2)}(\Q_p)\cap \overline{J(\Q)}$. Their existence is guaranteed when the analogous Chabauty assumption 
\begin{align}
\label{chabautycondition}
r<g-1
\end{align}
is satisfied. (In general, for $d$-fold symmetric powers this is $r<g-(d-1)$.)
However, unlike in the classical case for $X^{(1)}$, finding two linearly independent vanishing differentials alone need not yield an effective upper bound for the number of rational points on $X^{(2)}$. In fact, even when (\ref{chabautycondition}) is satisfied, $X^{(2)}(\Q)$ may still be infinite due to the existence of a degree 2 map $X\to C$ to a curve $C$ with infinitely many rational points. This happens for $X=X_0(N)$ with $N\in \{43,53,61,65\}$. Moreover, even if $C(\Q)$ and $X^{(2)}(\Q)$ are both finite, Chabauty's method still fails to find  an upper bound for $X^{(2)}(\Q)$ when no vanishing differential exists on $C$. This happens for $X=X_0(N)$ with $N\in \{67,73\}$. For $N=57$, this symmetric Chabauty method does succeed in determining all quadratic points. 
\subsubsection{Precise formulation}
\label{precisesymm}
Let $\mathcal{X}$ be a proper $\Z_p$-scheme with generic fibre $X$, which is smooth over $\spec (\Z_p$). Smoothness here simply means that $\widetilde{X}:=\mathcal{X}_{\F_p}$ is also non-singular, c.f. Proposition 2.9 in \cite{advancedsilverman}. This is in fact the minimal proper regular model of $X$; see the section on schemes in \cite{hindry}. Note that in particular we assume that $X$ has good reduction at $p$. 

 Inside the global differential forms $\Omega_{X/\Q_p}(X)$, we have the sub-$\Z_p$-module $\Omega_{\mathcal{X}/\Z_p}(\mathcal{X})$. Coleman integration defines a bilinear pairing
\begin{align}
\label{pairing}
\Omega_{X/\Q_p}(X)\times J(\Q_p) \to \Q_p, \;\; (\omega,\left[\sum_i P_i-Q_i\right])\mapsto \sum_i\int_{Q_i}^{P_i}\omega,
\end{align}
with kernel  equal to $J(\Q_p)_{\mathrm{tors}}$ on the right and 0 on the left. Define $V$ to be the annihilator of $J(\Q)$ with respect to the pairing (\ref{pairing}), and $\mathcal{V}:=V\cap \Omega_{\mathcal{X}/\Z_p}(\mathcal{X})$.  Let $\widetilde{V}$  be the image of $\mathcal{V}$  under the reduction map $\Omega_{\mathcal{X}/\Z_p}(\mathcal{X})\to \Omega_{\widetilde{X}/\F_p}(\widetilde{X})$ (see the proof of Lemma \ref{surjectivelemma} for a definition of this reduction map).

A priori $\widetilde{V}$ is, despite being a space of differentials on the reduction $\widetilde{X}$, defined in terms of the differentials on $\mathcal{X}$. In Section \ref{vanishingdiffs} we show that, when $X=X_0(N)$ for $N$ in our list, we can give an alternative description of $\widetilde{V}$ in terms of data on the reduction $\widetilde{X}$ only. This saves us from having to compute expansions of the vanishing differentials on $\mathcal{X}$, which involves the computationally non-trivial problem of lifting uniformisers around the point of expansion (see Section 2.5 in \cite{siksek}). 

Consider $\mathcal{Q}=\{Q_1,Q_2\}\in X^{(2)}(\Q)$. If $Q_1=Q_2$, we say that $Q_1$ has multiplicity $m=2$; else $Q_1$ has multiplicity $m=1$. Now choose a basis $\omega_1,\ldots, \omega_k$ for $\widetilde{V}$. Also choose a prime $v$ of $\Q(Q_1)$ above $p$, reductions respect to which we also denote by a tilde. Choose uniformisers $t_{\widetilde{Q}_j}$ at $\widetilde{Q}_j$ for each $j$. Now for each $i$ and $j$, we can expand $\omega_i$ around $\widetilde{Q}_j$ as a formal power series:
\[
\omega_i=(a_0(\omega_i,t_{\widetilde{Q}_j})+a_1(\omega_i,t_{\widetilde{Q}_j})t_{\widetilde{Q}_j}^1+\ldots)\dd t_{\widetilde{Q}_j}.
\]
When $m=1$, define 
\[
\widetilde{\mathcal{A}}:=\begin{pmatrix}
a_0(\omega_1,t_{\widetilde{Q_1}}) & a_0(\omega_1,t_{\widetilde{Q_2}}) \\
\vdots & \vdots \\
a_0(\omega_k,t_{\widetilde{Q_1}}) & a_0(\omega_k,t_{\widetilde{Q_2}}) 
\end{pmatrix}
\]
and when $m=2$ and $p>2$, set
\[
\widetilde{\mathcal{A}}:=\begin{pmatrix}
a_0(\omega_1,t_{\widetilde{Q_1}}) & a_1(\omega_1,t_{\widetilde{Q_1}})/2 \\
\vdots & \vdots \\
a_0(\omega_k,t_{\widetilde{Q_1}}) & a_1(\omega_k,t_{\widetilde{Q_1}})/2
\end{pmatrix}.
\]

\begin{theorem}[Symmetric Chabauty, Siksek]
 Suppose that $p>2$ and that $p\neq 3$ when $[\F_p(\widetilde{Q_1}):\F_p]=1$. 
If $\mathrm{rank} (\widetilde{\mathcal{A}})=2$ then $\mathcal{Q}$ is the unique point of $X^{(2)}(\Q)$ in its residue class modulo $p$.
\label{theorem1}
\end{theorem}
\begin{proof}
This is a special case of Theorem 1 in \cite{siksek} except for one difference: we work with expansions of differentials on the reduction $\widetilde{X}$, rather than reducing the coefficients of the expansion of differentials on $\mathcal{X}$.  To justify this, note that expansions of the vanishing differentials on $\mathcal{X}$ in \cite{siksek} are taken with respect to uniformisers $t_{Q_j}$ at the points $Q_j$ that also reduce to uniformisers $t_{\widetilde{Q}_j}$ at the reduced points. So if $\eta \in \Omega_{\mathcal{X}/\Z_p}(\mathcal{X})$ has expansion
\[
\eta=(a_0+a_1t_{Q_j}+\ldots)\dd t_{Q_j}
\]
then the expansion of $\widetilde{\eta}$ with respect to $t_{\widetilde{Q}_j}$ has the reductions $\widetilde{a}_i$ as coefficients. Therefore,  the reduction of the matrix $\mathcal{A}$ in \cite{siksek} equals our $\widetilde{\mathcal{A}}$, at least up to a change of uniformiser and possibly the removal of redundant rows (in case some vanishing differentials agree modulo $p$, c.f. Remark \ref{p=2remark}). These changes leave the rank unaffected. 
\end{proof}
Maarten Derickx \cite{dkss} has written down a similar criterion in terms of the differentials on $\widetilde{X}$ to show points are alone in their residue classes. His approach using the theory of formal immersions was already implicit in the works of Mazur \cite{mazur5} and Kamienny \cite{kamienny2}. It has the advantage of not introducing denominators in the matrix $\widetilde{\mathcal{A}}$ and can thus also be applied to $p=2$ when $m=2$. On the other hand, Siksek's method allows slightly more flexibility in the choice of the vanishing differentials.

\begin{remark} It would be interesting to determine if the combination of sufficiently many vanishing differentials for $X^{(2)}$ as well as for the curves $C$ admitting a degree 2 map $X\to C$ suffices for the Chabauty method to succeed on $X^{(2)}$, and whether a more intrinsic condition on $X$ exists; but we do not pursue this here. 
\end{remark}

\subsection{The relative case}\label{relativechab}
\subsubsection{Introduction}
When $X=X_0(N)$, we found that, in the cases where Chabauty fails, it is indeed due to a degree 2 map $X\to C$ for some curve $C$. Moreover, except for $N=37$ (a case treated separately in Section \ref{hypersection}), there is a single such curve $C$. We thus slightly change perspective and use a relative version of Chabauty's method to instead determine all \emph{exceptional points} of $X^{(2)}(\Q)$, i.e. those that do not come from $C(\Q)$ via the degree 2 map $X\to C$. In each case $C(\Q)$ can be described separately, either because it is an elliptic curve and we can find generators (this happens for $X_0(N)$ with $N\in \{43,53,61,65\})$ or because it is a hyperelliptic curve and its rational points have already been determined using quadratic Chabauty in \cite{bbbmtv} (this happens for $X_0(N)$ with $N\in \{67,73\}$). 

The main idea behind the relative Chabauty method is the following. Suppose that we have a degree 2 map $X\to C$, where $C$ is a curve with good reduction at the prime $p$. Then $J$ is isogenous to $J(C)\times A$ for some Abelian variety $A$, and we have a trace map $\mathrm{Tr}: \Omega_{X/\Q_p}(X)\to \Omega_{C/\Q_p}(C)$. Denote the map $J\to A$ by $\pi$. As long as we can find sufficiently many vanishing differentials $\omega$ whose image (also denoted by $\omega$) in $\Omega_{J/\Q_p}(J)$ is in the subspace $\pi^*\Omega_{A/\Q_p}(A)$,  Chabauty's method is expected to succeed in determining the exceptional points using only these differentials. In terms of $X$, the condition $\omega \in \pi^*\Omega_{A/\Q_p}$ means that $\omega$ is in the kernel of the trace map. In practice $A(\Q)$ will have rank zero and all differentials in $\pi^*\Omega_{A/\Q_p}(A)$ are vanishing differentials, c.f. Lemma \ref{vanishingdiffs}.

\begin{remark}\label{rationalrk}
Suppose that we can find one rational point $P$ on $X$. Then each other rational point $Q\in X(\Q)$ will occur in $X^{(2)}(\Q)$ at least twice: as $\{Q,Q\}$ and as $\{P,Q\}$. At most one of these pairs is the pullback of the image of $Q$ under $X\to C$, so at least one pair will occur as an exceptional point of $X^{(2)}(\Q)$. This means that, as a byproduct, we also determine the rational points when we find these exceptional points. Indeed each of the modular curves we consider contains a rational point.
\end{remark}

\subsubsection{Precise formulation}
We continue with the notation from Section \ref{precisesymm}. Assume that $X$ admits  a degree 2 map $\rho: X\to C$ to a curve $C$ and consider a proper smooth $\Z_p$-scheme $\mathcal{C}$ with $C$ as generic fibre. We assume that $\rho$ extends to a morphism $\mathcal{X}\to \mathcal{C}$, thus obtaining a corresponding reduced degree 2 map $\widetilde{X}=\mathcal{X}_{\F_p}\to \mathcal{C}_{\F_p}$ between non-singular curves. Note that in particular we assume that $C$ has good reduction at $p$. We denote the points of $X^{(2)}(\Q)$ coming from rational points of $C$ by $\rho^*C(\Q)$. 

Define $\mathcal{V}_0:=\mathcal{V}\cap \mathrm{ker}(\mathrm{Tr}: \Omega_{X/\Q_p}(X)\to \Omega_{C/\Q_p}(C)$ and let $\widetilde{V}_0$ be its image under the reduction map. As in the non-relative case, we show in Section \ref{vanishingdiffs} how $\widetilde{V}_0$ can be computed in terms of $\widetilde{X}$ only. 

We now consider a point $\mathcal{Q}=\{Q_1,Q_2\}\in \rho^*C(\Q)$ and choose a basis $\omega_1,\ldots,\omega_k$ for $\widetilde{V}_0$. Define $N$ as in Section \ref{precisesymm} and again let $t_{\widetilde{Q}_1}$ be a uniformiser at the reduction $\widetilde{Q}_1$ of $Q_1$ at the chosen prime $v$ above $p$.

\begin{theorem}[Relative symmetric Chabauty, Siksek]
\label{theorem2}
Suppose that $p>2$ when $[\F_p(\widetilde{Q_1}):\F_p]=1$.
If there exists $i\in \{1,\ldots,k\}$ such that
\[
\frac{\omega_i}{\dd t_{\widetilde{Q}_1}}|_{t_{\widetilde{Q}_1}=0} \neq 0
\]
then every point in $X^{(2)}(\Q)$ belonging to the residue class of $\mathcal{Q}$ comes from a point in $C(\Q)$. 
\end{theorem}
\begin{proof}
This is a special case of Theorem 2 in \cite{siksek}, with the same difference as in Theorem \ref{theorem1}; see its proof for the justification. 
\end{proof}
Note that $\frac{\omega_i}{\dd t_{\widetilde{Q}_1}}|_{t_{\widetilde{Q}_1}=0}$ is simply the constant coefficient of the expansion of $\omega_i$ at $t_{\widetilde{Q}_1}$, and the non-zero condition says that the corresponding $k\times 1$ matrix has rank 1. 

It is the combination of Theorem \ref{theorem1} for the points in $X^{(2)}(\Q)\setminus \rho^*C(\Q)$ and Theorem \ref{theorem2} for the points in $\rho^*C(\Q)$ that provably determines $X^{(2)}(\Q)$ for all curves $X$ we consider.

\subsection{The Mordell--Weil sieve}\label{sieve}
Again consider a (smooth) non-singular curve $X/\Q$ with Jacobian $J$. If $X$ admits a degree 2 map to another curve over $\Q$, we denote it by $\rho: X\to C$. 

In Chabauty's method the choice of a prime $p$ (of good reduction for $X$ and, if appropriate, also for $C$) is free to choose. The happy outcome is a number of $p$-adic discs in $X^{(2)}$ where no unknown rational points can lie. The Mordell--Weil sieve is a way of combining this $p$-adic information for several primes $p$ to make sure all of $X^{(2)}$ is covered. We have implemented the sieve described in Section 5 of \cite{siksek}, with one minor difference to be pointed out in Remark \ref{differenceremark}. 

Choose primes $p_1,\ldots,p_r$ (see the next Section on how to choose wisely). We assume given the following input:
\begin{itemize}
\item[(i)] a number of divisors $D_1,\ldots, D_n$ generating a subgroup $G$ of $J(\Q)$ of finite index,
\item[(ii)] a positive integer $I$ such that $I\cdot J(\Q)\subset G$, and
\item[(ii)] a (finite) list $\mathcal{L'}$ of known rational points for $X^{(2)}$, which, in the case of a degree 2 map $X\to C$, may also include points of $\rho^*C(\Q)$. 
\end{itemize}
In case of a degree 2 map $X\to C$, define $\mathcal{L}:=\mathcal{L}'\cup \rho^*C(\Q)$. Else, let $\mathcal{L}:=\mathcal{L}'$.

Suppose that $\mathcal{P}$ is a hypothetical point in $X^{(2)}(\Q)\setminus \mathcal{L}$. Our objective is to show that such $\mathcal{P}$ cannot exist. Let $\phi: \Z^n\to J(\Q)$ be the map given by $\phi(a_1,\ldots,a_n)=\sum_i a_i D_i$. It has image $G$. We choose a rational degree 2 divisor $\infty$  and define $\iota: X^{(2)}(\Q)\to G$ by $\mathcal{Q}\mapsto I\cdot [\mathcal{Q}- \infty]$. Let $p$ be one of the primes $p_i$. We denote reduction modulo $p$ by a tilde. Also define $\iota_p: X^{(2)}(\F_p)\to J(\F_p)$ by $\mathcal{R}\mapsto I\cdot [\mathcal{R}-\widetilde{\infty}]$. Finally define $\phi_p$ to make the diagram
\[
\xymatrixcolsep{6pc}\xymatrix{
\mathcal{L} \ar[r]^{i} \ar[rd] & X^{(2)}(\Q) \ar[r]^{\iota} \ar[d]^{\mathrm{red}} & G \ar[d]^{\mathrm{red}} & \ar[l]_{\phi} \ar[ld]^{\phi_p} \Z^n \\ & \widetilde{X}^{(2)}(\F_p) 
  \ar[r]_{\iota_p}& J(\F_p) & 
}
\]
commute.  
Note that $\iota_p(\widetilde{\mathcal{P}})=\mathrm{red}(\iota(\mathcal{P}))$, hence $\iota_p(\widetilde{\mathcal{P}})\in \mathrm{Im}(\phi_p)$. Since $\mathcal{P}\notin \mathcal{L}$, our Chabauty method limits the possible values of $\widetilde{\mathcal{P}}$, which in turn reduces the union $W_p$ of $\ker(\phi_p)$-cosets of $\Z^n$ that could possibly be mapped to $\iota_p(\widetilde{\mathcal{P}})$ under $\phi_p$. The key observation here is that $\Z^n$ is independent of $p$, and that these $W_p$ can thus be compared for different values of $p$. In particular, if they have empty common intersection, then our hypothetical point $\mathcal{P}$ cannot exist.

\begin{definition}
\label{sievedef}
Define $\mathcal{M}_p\subset \widetilde{X}^{(2)}(\F_p)$ to be the set of points  $\mathcal{R}\in \iota_p^{-1}(\mathrm{Im}(\phi_p))$ satisfying one of the following:
\begin{itemize}
    \item[(i)] $\mathcal{R}\notin \mathrm{red}(\mathcal{L'})$,
    \item[(ii)] $\mathcal{R}=\widetilde{\mathcal{Q}}$ for some $\mathcal{Q}\in \mathcal{L'}$ \emph{not} satisfying the conditions of Theorem \ref{theorem1} or
    \item[(iii)] there is a degree 2 map $X\to C$ and $\mathcal{R}=\widetilde{\mathcal{Q}}$ for some $\mathcal{Q}\in \mathcal{L'}\cap \rho^*C(\Q)$ \emph{not} satisfying the conditions of Theorem \ref{theorem2}.
\end{itemize}
\end{definition}
Then by construction, the reduction $\widetilde{\mathcal{P}}$ of our hypothetical point is in $\mathcal{M}_p$. We conclude the following.
\begin{theorem}[Mordell--Weil sieve]
If 
\[
\bigcap_{i=1}^r\phi_{p_i}^{-1}(\iota_{p_i}(\mathcal{M}_{p_i}))=\emptyset
\]
then $X^{(2)}(\Q)=\mathcal{L}$.
\end{theorem}
\begin{remark}
For each prime $p$, we note that $\widetilde{X}^{(2)}(\F_p)$, $J(\widetilde{X})(\F_p)$, $\iota_p$, $\phi_p$ and $\mathrm{red}|_{\mathcal{L}}$ can all be computed explicitly. Therefore, each coset $\phi_p^{-1}(\iota_p(\mathcal{M}_{p_i}))$ can be determined explicitly. The computations we perform with cosets of subgroups of $\Z^n$ require only linear algebra over $\Z$ and can be done almost instantaneously by computer algebra systems.
\end{remark}
\begin{remark}
\label{differenceremark}
The only difference between this sieve and the one described in Section 5 of \cite{siksek} is that we work with cosets in $\Z^n$ rather than cosets of increasingly large finite quotients of $\Z^n$. This saves us from an explosion of the size of such finite quotients caused by the Chinese Remainder Theorem.
\end{remark}

\subsubsection{Prime-choosing heuristics}
An interesting aspect of the Mordell--Weil sieve is the question of which primes to choose in which order. We have found the naive choice of a small number of small primes to be sufficient when $N\notin \{67,73\}$. For the remaining two cases, however, we used the prime-choosing heuristics as described in this section. These were inspired by those in Section 11 of \cite{bugeaud}. Note that, even though this strategy works for us, it may not be optimal. For example, we attempt to choose the primes one by one in a near-optimal sense, rather than attempting to find a near-optimal \emph{set} of primes. A more detailed discussion about choosing primes in the Mordell--Weil sieve can be found in \cite{bruinstoll}.

Suppose that we have already chosen the primes $p_1,\ldots,p_k$ and we would like to choose the next prime $p_{k+1}$. Let us write $W_k:=\bigcap_{i=1}^k\phi_{p_i}^{-1}(\iota_{p_i}(\mathcal{M}_{p_i}))$. This is a union of $A_k$-cosets, where $A_k\subset \Z^n$ is the subgroup $\bigcap_{i=1}^k \ker \phi_{p_i}\subset \Z^n$.

\begin{choice}
We first establish a bounded range of primes $\mathbb{P}$ in which computations are reasonably fast, e.g. all primes up to 100. Then we choose $p_{k+1}$ to be the prime of good reduction in $\mathbb{P}\setminus \{p_1,\ldots,p_k\}$ minimising 
\[
\frac{[A_k: \ker (\phi_{p_{k+1}}|_{A_k})]}{p_{k+1}^{g-2}}.
\]
\end{choice}

Let us explain the ideas behind this choice. Let $w_k$ be the number of $A_k$-cosets in $W_k$. We aim to choose the next prime $p_{k+1}$ in such a way that $w_{k+1}$ is as small as possible. Note that, after choosing $p_{k+1}$, $W_{k+1}=W_k\cap \phi_{p_{k+1}}^{-1}(\iota_{p_{k+1}}(\mathcal{M}_{p_{k+1}}))$ consisting of $A_{k+1}=\ker (\phi_{p_{k+1}}|_{A_k})$-cosets. So a priori, $w_k$ gets multiplied by a factor 
\[
[A_k: \ker (\phi_{p_{k+1}}|_{A_k})].
\]
Next, to form $W_{k+1}$ we remove the $A_{k+1}$-cosets not mapping to $\iota_{p_{k+1}}(\widetilde{X}^{(2)}(\F_p))$. Since $\widetilde{X}^{(2)}(\F_p)\to J(\widetilde{X})(\F_p),\; \mathcal{Q}\mapsto [\mathcal{Q}-\infty]$ is injective, one would expect (assuming the image of $\widetilde{X}^{(2)}(\F_p)$ in $J(\widetilde{X})(\F_p)$ is randomly behaved) a proportion
\[
\frac{\#\widetilde{X}^{(2)}(\F_p)}{\#J(\widetilde{X})(\F_p)}
\]
of the cosets to remain.

Finally, we also remove the $A_{k+1}$-cosets mapping to elements $D\in J(\widetilde{X})(\F_p)$ such that each point in $\iota_{p_{k+1}}^{-1}(D)$  is excluded by Chabauty. This is a contribution that we decide to neglect when choosing our prime, both because it is hard to estimate and because we work with a fixed finite set $\mathcal{L}'\subset \mathcal{L}$. As the size of $p$ increases, $\mathrm{red}(\mathcal{L'})$ thus forms an increasingly small proportion of $\widetilde{X}^{(2)}(\F_p)$.  

Next, we note that $\#J(\F_p)=p^g+o(p^g)$ as $p\to \infty$ and $\widetilde{X}^{(2)}(\F_p)=p^2+o(p^2)$ as $p\to \infty$ by classical bounds. Hence we approximate $\#\widetilde{X}^{(2)}(\F_p)/\#J(\F_p)$ by $p^{2-g}$. Therefore, we approximate $w_{k+1}$ by
\[
\frac{[A_k: \ker (\phi_{p_{k+1}}|_{A_k})]}{p_{k+1}^{g-2}}w_k,
\]
which explains our choice.

Note that we do need to make a little precomputation before choosing each new prime, but this precomputation is relatively fast compared to checking the Chabauty method at each point. Moreover, the relatively expensive part of this precomputation is computing $J(\F_p)$ and $\phi_p$ for each $p\in \mathbb{P}$, which only needs to be done once.

\section{The input}
We continue with the notation from the previous Section. When $X=X_0(N)$ for $N \in \{43, 53, 61, 57, 65, 67,73\}$, we always have a degree 2 map $X_0(N)\to X_0(N)^+$, where $X_0(N)^+$ is the quotient of $X_0(N)$ by the Atkin--Lehner involution $w_N$. Except for $N=57$, we will use relative symmetric Chabauty with respect to $C=X_0(N)^+$. In this Section we describe how we computed the necessary input to make this work. This consists of the following four parts:
\begin{enumerate}
\item[(i)] explicit defining equations for $X_0(N)$ as well as the corresponding Atkin--Lehner involutions,
    \item[(ii)] a list $\mathcal{L}'$ of known quadratic points on $X$, including points coming from $C(\Q)$ in the case of a degree 2 map $X\to C$,
    \item[(iii)] a list of generators of a subgroup $G\subset J(\Q)$ and an integer $I$ such that $I\cdot J(\Q)\subset G$,
    \item[(iv)] at least 2 linearly independent vanishing differentials for each prime $p$, which are in the kernel of $\mathrm{Tr}: \Omega_{X/\Q_p}(X)\to \Omega_{C/\Q_p}(C)$ in case of a degree 2 map $X\to C$.
\end{enumerate}

In this Section we describe how we obtain this input. We note that knowing more than two vanishing differentials will improve the speed and the chance of success of the sieve.
\subsection{A model for $X_0(N)$}
The Small Modular Curves package in \texttt{Magma} is a great tool for computing models of modular curves $X_0(N)$. As, however, $X_0(65),X_0(67)$ and $X_0(73)$ are not in this database, we instead compute models for all our modular curves using the code written by Ozman and Siksek; see Section 3 of \cite{ozman}. This, too, uses the canonical embedding to compute the models. Moreover, Ozman and Siksek also explicitly compute the action of the Atkin--Lehner operators on their models, as well as the equations for the $j$-invariant map $X_0(N)\to X(1)$.

For primes $p$ where the chosen model has good reduction, we define $\mathcal{X}_0(N)$ to be the $\Z_p$-scheme defined by these equations. In these cases, $\mathcal{X}_0(N)$ is our minimal proper regular model, c.f. Section 2.2.2. 
\subsection{Searching for quadratic points} 
We assume  henceforth that we have found a model for $X$ in $\P^k$, with coordinates $x_0,\ldots,x_k$. We search for quadratic points on $X$ by intersecting $X\subset \P^k$ with hyperplanes of the form 
\[
b_0x_0+\ldots + b_kx_k=0,
\]
where $b_0,\ldots,b_k\in \Z$ are chosen coprime and up to a certain bound. When the decomposition of the divisor (over $\Q$) corresponding to this intersection contains effective degree 2 divisors, we have found quadratic points. In practice it sufficed for us to consider only the hyperplanes defined by $|b_i|\leq 10$ for $i\in \{1,\ldots,k\}$. Unlike searching for rational points, this can be a time consuming exercise: for the genus 5 curve $X_0(67)\subset \P^4$, for example, the search took multiple hours. This is largely due to the time taken to decompose these hyperplane intersections into linear combinations of irreducible effective divisors.

\subsection{Determining subgroups of the Mordell--Weil groups}
We first consider the general case, where $X/\Q$ is again any (projective non-singular) curve, $\Gamma$ is a finite subgroup of $\mathrm{Aut}_{\Q}(X)$ and $C=X/\Gamma$. Denote the quotient map by $\rho: X\to C$. This map has degree $\deg\rho = \#\Gamma$. We still assume that $X$ has a rational point. We now denote the Jacobian of $X$ by $J(X)$ to emphasize the difference with the Jacobian $J(C)$ of $C$. After choosing compatible base points for the maps $\iota_X: X\to J(X)$ and $\iota_C: C\to J(C)$, we then obtain a commuting diagram
\sqcommdiag{X}{J(X)}{C}{J(C).}{\rho}{\iota_X}{\rho_*}{\iota_C}
%

We now make the extra assumption that
$\mathrm{rk} J(X)(\Q)=\mathrm{rk} J(C)(\Q)$; denote this common rank by $r$. Let $G$ be the subgroup of $J(X)(\Q)$ generated by $J(X)(\Q)_{\mathrm{tors}}$ and $\rho^*J(C)(\Q)$. 

\begin{proposition}
We have that $(\#\Gamma)\cdot J(X)(\Q)\subset G$.
\label{mwgroup}
\end{proposition}
\begin{proof}
Let $P_1,\ldots,P_r \in J(C)(\Q)$ be linearly independent generators  of $J(C)(\Q)/J(C)(\Q)_{\mathrm{tors}}$. Set $D_i=\rho^*P_i$ for each $i$. Then $\rho_*D_i=\deg(\rho)P_i=(\#\Gamma)P_i$ and hence $D_1,\ldots, D_r$ are also linearly independent.

We consider each $D_i$ as an element of $J(X)(\Q)/J(X)(\Q)_{\mathrm{tors}}$ and each $P_i$ as an element of  $J(C)(\Q)/J(C)(\Q)_{\mathrm{tors}}$, while still denoting the maps between these quotients by $\rho_*$ and $\rho^*$. As $\rho_*D_i=(\#\Gamma)P_i$ for each $i$, we see that $\rho_*G=(\#\Gamma)\cdot J(C)(\Q)/J(C)(\Q)_{\mathrm{tors}}$. In particular, as we assume $\rho_*$ to be injective modulo torsion (equal rank), we find that $(\#\Gamma)\cdot J(X)(\Q)/J(X)(\Q)_{\mathrm{tors}}\subset G/J(X)(\Q)_{\mathrm{tors}}$. 
%
%
\end{proof}

We now consider the special case $X=X_0(N)$ for $N$ in our list. Let us first compute the torsion subgroup of the Jacobian. For each $N$, let $C_0(N)$ be the subgroup of $J_0(N)(\overline{\Q})$ generated by the classes of differences of cusps, and $C_0(N)(\Q)$ its subgroup fixed by $\mathrm{Gal}(\overline{\Q}/\Q)$. This is called the \emph{rational cuspidal subgroup}. By the Manin-Drinfeld theorem \cite{manin}, \cite{drinfeld}, $C_0(N)(\Q)\subset J_0(N)(\Q)_{\mathrm{tors}}$ for each $N$, and a conjecture of Ogg proved by Mazur \cite{mazur1} tells us that $J_0(N)(\Q)_{\mathrm{tors}}=C_0(N)(\Q)$ for prime values of $N$. Moreover, Mazur also showed in that case that the order of $J_0(N)(\Q)_{\mathrm{tors}}$ is the numerator of $(p-1)/12$. Recent results (\cite{oggconj1},\cite{oggconj2},\cite{oggconj3},\cite{oggconj4},\cite{ozman}) have verified the equality $C_0(N)(\Q)=J_0(N)(\Q)_{\mathrm{tors}}$ for various non-prime values of $N$, to which we can now add 57 and 65.  
%
\begin{lemma}
For $N\in \{43,53,61,67,73\}$ the torsion subgroup of $J_0(N)(\Q)$ is generated by the difference of the two cusps, the orders of which are, respectively, 7, 13, 5, 11 and 6. The torsion subgroups of $J_0(57)(\Q)$ and $J_0(65)(\Q)$ also equal their rational cuspidal subgroups, which are $\Z/6\Z\times \Z/30\Z$ and $\Z/2\Z\times \Z/84\Z$ respectively as abstract groups.
\label{torsionlemma}
\end{lemma}
\begin{proof}
For the prime values of $N$, this is the aforementioned theorem of Mazur. For $N=57$, we use \texttt{Magma} code of Ozman and Siksek \cite{ozman} to compute $C_0(57)(\Q)$ as a group, which gives $\Z/6\Z\times \Z/30\Z$. As 5 is a prime of good reduction for $X_0(N)$, we find that $J_0(57)(\Q)_{\mathrm{tors}}$ injects into $J_0(57)(\F_5)$.
We compute that
\[
 J_0(57)(\F_5)\simeq \Z/3\Z\times \Z/6\Z\times \Z/180\Z
\]
so that the index $[J_0(57)(\Q)_{\mathrm{tors}}:C_0(57)(\Q)]$ divides 18. Looking at $J_0(57)(\F_{23})$, we find that the index is coprime to 3, and from $J_0(57)(\F_{11})$ we learn that the index is odd, leaving index 1 as only option. A similar computation at the primes 3 and 11 allows us to determine $J_0(65)(\Q)_{\mathrm{tors}}$. 
\end{proof}
%
For $N\in \{43,53,65,61,67,73\}$, let $\Gamma$ be the subgroup of $\mathrm{Aut}(X_0(N))$ generated by $w_N$, and for $N=57$ let $\Gamma=\langle w_{19},w_{3}\rangle$. Define $C(N):=X_0(N)/\Gamma$. For each $N\neq 57$, this curve will be the degree two quotient of $X_0(N)$ denoted before in greater generality by $C$. Note that $C(N)=X_0(N)^+$ for $N\neq 57$ and $C(N)$ is the curve often denoted by $X_0(N)^*$ for $N=57$. We verify that, with our choice of $X=X_0(N)$ and $\Gamma$, the hypothesis for Proposition \ref{mwgroup} holds true.
\begin{lemma}
\label{ranklemma}
For each $N$ in our list, we have 
$
\mathrm{rk}J_0(N)(\Q) = \mathrm{rk} J(C(N))(\Q).
$
\end{lemma}
\begin{proof}
We checked this using the algorithm of Stein \cite{stein} as implemented in the modular arithmetic geometry package in \texttt{Magma}. For each $N$, $J_0(N)\sim \prod_f A_f$ is isogenous to a product of simple modular abelian varieties $A_f$ corresponding to Galois orbits $f$ of newforms of weight 2 and level $N$. Now $J(C(N))$ is isogenous to the product of those $A_f$ where $f$ is invariant under $\Gamma$. For each $N$ we found $J(C(N))$ to be simple, hence $J(C(N))\sim A_f$ for some orbit $f$.  We checked that the values of the $L$-series $L(g,s)$ are non-zero at $s=1$ for the eigenforms $g$ not conjugate to $f$. This means that the corresponding $A_g$ have analytic rank 0, and hence algebraic rank zero by a theorem of Kolyvagin and Logachev \cite{kolyvagin}.
\end{proof}
Write $r(N):=\mathrm{rk} J_0(N)(\Q)$. We note that $r(N)=2$ for $N\in \{67,73\}$ and $r(N)=1$ for the other values of $N$. Also $C(N)$ has genus 2 for $N\in \{67,73\}$ and genus 1 otherwise, making $C(N)$ either an elliptic curve or a hyperelliptic curve. When $C(N)$ is an elliptic curve, we can compute a basis for its Mordell--Weil group using Cremona's algorithm \cite{cremona} implemented in \texttt{Magma}. When it is hyperelliptic of genus 2, we manage to do the same using height bounds and an algorithm of Stoll \cite{stoll}. We have written down the explicit generators in Section \ref{results}. 
%
%
Together with Lemma \ref{torsionlemma}, this gives us a complete description of a suitable subgroup of $J_0(N)(\Q)$ for each of the values of $N$ we consider.


%

\subsection{The vanishing differentials}\label{vanishingdiffs}
For $N\neq 57$, let $V$ be the image of $1-w_N^*: \Omega_{X_0(N)/\Q}(X_0(N))\to \Omega_{X_0(N)/\Q}(X_0(N))$. For $N=57$, let $V\subset \Omega_{X_0(57)/\Q}(X_0(57))$ be the sum of the images of $1-w_3^*$ and $1-w_{19}^*$. 
\begin{lemma}
All differentials $\omega\in V$ annihilate  $J_0(N)(\Q)$ via the integration pairing, i.e.
\[
\int_{0}^{D}\omega =0 \text{ for all }D\in J_0(N)(\Q).
\]
Moreover, $V$ is contained in the kernel of $\mathrm{Tr}: \Omega_{X_0(N)/\Q}(X_0(N))\to \Omega_{C(N)/\Q}(C(N))$. 
\end{lemma}
\begin{proof}
Note that, as the Atkin--Lehner maps are involutions, the image of $1-w_N^*$ is the kernel of $1+w_N^*$, which is the trace map for $N\neq 57$. For $N=57$, both the images of $1-w_3^*$ and $1-w_{19}^*$ are in the kernel of the trace map $1+w_3^*+w_{19}^*+w_{57}^*$. If we identify $V$ with a subspace of $\Omega_{J_0(N)/\Q}(J_0(N))$, this means that $V\subset \pi^*\Omega_{A/\Q}$, where $J_0(N)\sim J(C)\times A$ as in the proof of Lemma \ref{ranklemma} and $\pi$ is the map $J_0(N)\to A$. It can also be seen directly using the algorithm of Stein \cite{stein} implemented in the modular abelian varieties package in \texttt{Magma} that the image of $1-w_N^*: J_0(N)\to J_0(N)$ is $A$ for $N\neq 57$, and that the union of the images of $1-w_3^*$ and $1-w_{19}^*$ is $A$ for $N=57$. Let $\omega\in V$ and consider $\eta\in \Omega_{A/\Q}$ such that $\omega=\pi^*\eta$. By Lemma \ref{ranklemma}, we find that $A(\Q)$ is torsion. Consider $D\in J_0(N)(\Q)$ and let $n$ be the order of $\pi(D)$ in $A(\Q)$. By the additive property of Coleman integration, we find
\[
\int_0^{D}\omega=\int_0^{\pi(D)}\eta =\frac{1}{n}\int_0^{n\pi(D)}\eta =0,
\]
as desired.
\end{proof}

As we have computed equations for $X_0(N)$ as well as the Atkin--Lehner involutions, we can compute $V$ explicitly for each $N$. Note that $\dim(V)=\mathrm{genus}(X)-\mathrm{genus}(C)\geq 2$ for every value of $N$ we consider. This indicates that the relative symmetric Chabauty method may succeed.

However, in order to apply Theorems \ref{theorem1} and \ref{theorem2}, we need to compute the image $\widetilde{V}$ of $\mathcal{V}=V\cap \Omega_{\mathcal{X}_0(N)/\Z_p}(\mathcal{X}_0(N))$ under the reduction map. To this end, consider a prime $p$ of good reduction for $X=X_0(N)$ and, when appropriate, for $C(N)$. Note that each Atkin--Lehner involution $w_M$ for $M\mid N$ is integral and thus extends to a map on the proper minimal regular model $\mathcal{X}_0(N)\to \spec \Z_p$. In particular, we obtain reduced maps $\widetilde{w}_M: \widetilde{X}_0(N)\to \widetilde{X}_0(N)$. The following proposition allows us to easily compute $\widetilde{V}$ in \texttt{Magma}.

\begin{proposition}
\label{vanishingprop}
For $N\neq 57$, $\widetilde{V}$ is the image of $1-\widetilde{w}_N^*: \Omega_{\widetilde{X}_0(N)/\F_p}(\widetilde{X}_0(N))\to \Omega_{\widetilde{X}_0(N)/\F_p}(\widetilde{X}_0(N))$. For $N=57$, $\widetilde{V}$ is the sum of the images of $1-\widetilde{w}_3^*$ and $1-\widetilde{w}_{19}^*$. 
\end{proposition}

The proof of this proposition relies on the following lemma, which in turn depends heavily on the fact that $X$ has good reduction at $p$.  This lemma is commonly known by experts, but we could not find a reference for it.

\begin{lemma}
\label{surjectivelemma}
Let $X/\Q$ be a (smooth projective) curve with minimal proper regular model $\mathcal{X}/\Z_p$. Then the reduction map $\Omega_{\mathcal{X}/\Z_p}(\mathcal{X})\to \Omega_{\widetilde{X}/\F_p}(\widetilde{X})$ is surjective.
\end{lemma}
\begin{proof}
 On the level of sheaves, the Cartesian diagram defining the fibre $\widetilde{X}$ yields the isomorphism
\[
\Omega_{\widetilde{X}/\F_p}\simeq \Omega_{\mathcal{X}/\Z_p}|_{\widetilde{X}}\otimes_{\Z_p}\F_p,
\]
where $\Omega_{\mathcal{X}/\Z_p}|_{\widetilde{X}}$ is the sheaf associated to 
\[
U\mapsto \lim_{\substack{\longrightarrow\\V\supset U}}\Omega_{\mathcal{X}/\Z_p}(V)
\]
(see Proposition 8.10 in \cite{hartshorne}). Moreover, this isomorphism defines the reduction map $\Omega_{\mathcal{X}/\Z_{p}}(\mathcal{X})\to \Omega_{\widetilde{X}/\F_p}(\widetilde{X})$ on global sections.

We first note that every Zariski-open subset of $\mathcal{X}$ containing $\widetilde{X}$ equals $\mathcal{X}$: if $V\subset \mathcal{X}$ is a Zariski-open set containing $\widetilde{X}$, then its complement is a closed subset (in $\mathcal{X}$) of the generic fibre $X$. However, the closure (in $\mathcal{X}$) of any point $P\in X$ contains its reduction $\widetilde{P}\in \mathcal{X}$, so $\mathcal{X}\setminus V=\emptyset$.

Now for simplicity, we write $A=\Omega_{\mathcal{X}/\Z_{p}}(\mathcal{X})$ and $B=\Omega_{\widetilde{X}/\F_p}(\widetilde{X})$. We will show that $A/pA\simeq B$. Note that $pA$ is in the kernel of $A\to B$. Next, by flatness of $\mathcal{X}\to \spec \Z_p$, we have that $\mathrm{rank}_{\Z_p} A=g =\dim_{\F_p}B$. Therefore, the induced map $A/pA\to B$ is a linear map of $\F_p$-vector spaces of equal dimension.  It thus suffices to show that each $\omega\in \ker (A\to B)$ is a $p$-fold. 
We first show this to be true locally. Consider $\omega\in \ker(A\to B)$. By the above isomorphism of sheaves, there exists an open cover $\{U\}$ of $\widetilde{X}$ such that 
\[
\Omega_{\widetilde{X}/\F_p}(U)\simeq\lim_{\substack{\longrightarrow\\V\supset U}}\Omega_{\mathcal{X}/\Z_p}(V)\otimes_{\Z_p}\F_p \;\;\; \text{ for each } U.
\]
Note that it does not matter whether we tensor with $\F_p$ inside or outside the limit. By this isomorphism, there must be open $V_U\supset U$ and $\eta_U\in \Omega_{\mathcal{X}/\Z_p}(V_U)$ such that $\omega|_{V_U}=p\eta_U$. Since $\bigcup_UV_U$ is an open subset of $\mathcal{X}$ containing $\widetilde{X}$, the $V_U$ must cover $\mathcal{X}$ entirely. As $(\omega|_{V_U})_{V_U}$ lifts to the global section $\omega$ and we are working in characteristic zero, also the $\eta_U$ must agree on overlaps and lift to $\eta\in A$ such that $\omega=p\eta$.
\end{proof}
\begin{proofofprop}
LLet $\widetilde{W}$ be the image of $1-\widetilde{w}_N^*$ for $N\neq 57$ and the sum of the images of $1-\widetilde{w}_3^*$ and $1-\widetilde{w}_{19}^*$ for $N=57$. We want to show that $\widetilde{W}=\widetilde{V}$. Note first that $\widetilde{V}\subset \widetilde{W}$ because if $\omega\in \ker(1-w_M^*)$ then $\widetilde{\omega}\in \ker (1-\widetilde{w}_M^*)$, so we need to prove the opposite inclusion. To this end, consider $\widetilde{\omega}\in \widetilde{W}$ and $\widetilde{\eta}$ such that $\widetilde{\omega}=(1-w_M^*)(\widetilde{\eta}$). Here $M=N$ for $N\neq 57$ and for $N=57$ we may assume $\widetilde{\omega}$ to be of this form for $M\in \{3,19\}$.  By the previous lemma, $\widetilde{\eta}$  lifts to some $\eta \in \Omega_{\mathcal{X}/\Z_p}(\mathcal{X})$. Then $(1-w_M^*)(\eta) \in \mathcal{V}$ reduces to $\widetilde{\omega}$, as desired.
%
\end{proofofprop}
\begin{remark}
\label{p=2remark}
When $p\neq 2$, the rank of $\mathcal{V}$ equals the dimension of $\widetilde{V}$: for such $p$, $\mathcal{V}$ is the set of those $\omega\in \Omega_{\mathcal{X}/\Z_p}(\mathcal{X})$ fixed by $\frac12(1-w_N^*)$, and if $\frac12(1-w_N^*)$ fixes $p\omega$ then it fixes $\omega$ too.  For $p=2$, however, $\frac12(1-w_N^*)$ may not map $\Omega_{\mathcal{X}/\Z_p}(\mathcal{X})$ into itself. In fact we found for $X=X_0(53)$ that $\mathrm{rank} (\mathcal{V})=3$  but $\dim \widetilde{V}=1$. This means that the matrix $\widetilde{\mathcal{A}}$ of Theorem \ref{theorem1} or the corresponding $1\times k$ matrix of Theorem \ref{theorem2} is reduced in size compared to the matrix $\mathcal{A}$ defined in Theorems 1 and 2 of \cite{siksek}, but note that this does not affect the rank.
\end{remark}
\section{Results for non-hyperelliptic $X_0(N)$}\label{results}
In this section we list all the (exceptional) quadratic points found for $X_0(N)$ with 
\[
N\in \{43,53,61,57,65,67,73\}
\]
as well as the elliptic curves they correspond to. We also list any non-cuspidal rational points, c.f. Remark \ref{rationalrk}. Models for $X_0(N)$ are always given in projective space. With $C$ we always mean the quotient of $X_0(N)$ by one or more Atkin--Lehner involutions such that $\mathrm{rk}J(C)(\Q) = \mathrm{rk}J_0(N)(\Q)$, and we denote the map $X\to C$ by $\rho$. For elliptic curves, $O$ always denotes the zero element, which equals $(0:1:0)$ in every case. The column denoted by CM lists the discriminant of the order by which the elliptic curve has complex multiplication if it does, and NO if it does not have CM. The column $\Q$-curve denotes the Atkin--Lehner operators defining an isogeny between the elliptic curve and its conjugate if there is one, and NO if it is not a $\Q$-curve.
\subsection{$X_0(43)$}
Model for $X_0(43)$:
\begin{align*}
x_0^3x_2 - 2x_0^2x_1^2 &+ 2x_0^2x_1x_2 - 2x_0^2x_2^2 + x_0x_1^3 + 3x_0x_1^2x_2 - 5x_0x_1x_2^2 \\&+ 3x_0x_2^3 - 9x_1^4 + 24x_1^3x_2 - 28x_1^2x_2^2 + 
    16x_1x_2^3 - 4x_2^4=0.
\end{align*}
Genus $X_0(43)$: 3.\\
Cusps: $(1:0:0)$, $(1:1:1)$.\\
$C=X_0(43)^+$: elliptic curve $y^2 + y = x^3 + x^2$ of conductor 43.\\
Group structure of $J(C)(\Q)$: $J(C)(\Q)= \Z\cdot [Q_C-O]$, where $Q_C:=(0 : -1 : 1)$.\\
Group structure of $G\subset J_0(43)(\Q)$: $G= \Z\cdot D_1\oplus \Z/7\Z\cdot D_{\mathrm{tor}}$, where 
$D_{\mathrm{tor}}=[(1:1:1)-(1:0:0)]$ and $D_1=[P+\overline{P}-(1:1:1)-(1:0:0)]=\rho^*[Q_C-O]$ for 
\[
P:=\left(\frac18(\sqrt{-7} + 3) : \frac18(-\sqrt{-7} + 5) : 1\right) \in X_0(43)(\Q(\sqrt{-7}))
\]
satisfying $\rho(P)=Q_C$.\\
Primes used in sieve: 5,7,11.\\
The following table contains a list of all quadratic points (up to Galois conjugacy) not coming from $X_0^+(43)(\Q)$ via $\rho: X_0(43)\to X_0^+(43)$ and all non-cuspidal rational points.
\begin{table}[h!]
    \centering
    \begin{tabular}{c c c c c c}
        Name & $\theta^2$ & Coordinates & $j$-invariant & CM & $\Q$-curve  \\ [0.5ex] 
        \hline \hline
        $P_0$ & - & $\left(0: \frac45 : \frac25 : 1\right)$ & -884736000 & -43 & -\\ 
        $P_1$ & -131 & $\left(\frac{1}{72}(\theta + 65) : \frac{1}{24}(\theta + 17) : 1\right)$ & \tiny{$\substack{(2646379314349235349704820159442238\theta \\+ 
    70713216722735823130811605070229181)\\ /3410605131648480892181396484375}$} & NO & NO \\   
    $P_2$ & -131 & $\left(\frac{1}{25}(-2\theta + 9) : 2/5 : 1\right)$ & \tiny{$\substack{(245508467396487686583118\theta \\- 
    5588464515419225929913951)\\/1641284836972685388135}$} &NO & NO \\   
    $P_3$ & -71 & $\left(\frac14(\theta + 1) : 1 : 0\right)$ & $\frac{-49\theta-977}{4}$ & NO & NO \\  
    $P_4$ & -71 & $ \left(\frac{1}{15}(-\theta + 8) : \frac{1}{30}(-\theta + 23) : 1\right)$ & \tiny{$\substack{(5111948195521623101849\theta\\ - 22519853936617719563123)\\/17592186044416}$} & NO & NO 
    \end{tabular}
    \label{table43}
\end{table}

\subsection{$X_0(53)$}
Model for $X_0(53)$: 
\begin{align*}
&9x _0^2  x_3 - 5  x_0  x_3^2 - 27  x_1^3 - 18  x_1^2  x_2 + 78  x_1^2  x_3 - 18  x_1  x_2^2 + 30 x_1  x_2  x_3 \\ &\quad\quad\quad\quad\quad\quad\;\;\;\;\;- 64  x_1  x_3^2 - 57  x_2^3 + 136  x_2^2  x_3 - 104  x_2  x_3^2 + 
    49  x_3^3=0,\\
    &x_0  x_2 - 2  x_0  x_3 - 3x_1^2 + 5  x_1  x_3 - 2 x_2^2 +   x_2  x_3 - 2  x_3^2=0.
\end{align*}
Genus $X_0(53)$: 4.\\
Cusps: $(1:0:0:0)$, $(1:1:1:1)$.\\
$C=X_0(53)^+$: elliptic curve $y^2 + xy + y = x^3 - x^2$ of conductor 53.\\\
Group structure of $J(C)(\Q)$: $J(C)(\Q)=\Z\cdot [Q_C-O]$, where $Q_C:=(0:-1:1)$.\\
Group structure of $G\subset J(X)(\Q)$: $G= \Z\cdot D_1\oplus \Z/13\Z\cdot D_{\mathrm{tor}}$, where $D_{\mathrm{tor}}=[(1:1:1:1)-(1:0:0:0)]$ and $D_1=[P+\overline{P}-(1:1:1:1)-(1:0:0:0)]=\rho^*[Q_C-O]$ for 
\[
P:=\left( 0, \frac16(-\sqrt{-11} + 5): 1: 1 \right)\in X_0(53)(\Q(\sqrt{-11}))
\]
satisfying $\rho(P)=Q_C$.\\
Primes used in sieve: 11,7. \\
There are {\bf no} quadratic points on $X_0(53)$ not coming from $X_0^+(53)(\Q)$ via $\rho: X_0(53)\to X_0^+(53)$ and {\bf no} non-cuspidal rational points. In particular, all quadratic points correspond to $\Q$-curves.

\subsection{$X_0(61)$}
Model for $X_0(61)$: 
\begin{align*}
& x_0^2  x_3 + x_0  x_1  x_3 - 2  x_0  x_3^2 - 2  x_1^3 - 6  x_1^2  x_2 + 5  x_1^2  x_3 - 5  x_1  x_2^2 + 4  x_1  x_2  x_3 - 6  x_2^3 + 14  x_2^2  x_3 - 11  x_2  x_3^2 + 4  x_3^3=0,\\
& x_0  x_2 - x_1^2 - x_1  x_2 - 2  x_2^2 + 2  x_2  x_3 - x_3^2=0
\end{align*}
Genus of $X_0(61)$: 4.\\
Cusps: $(1:0:0:0)$, $(1:0:1:1)$.\\
$C=X_0(61)^+$: elliptic curve $y^2 + xy = x^3 + 6x^2 + 11x + 6$ of conductor 61.\\
Group structure of $J(C)$: $J(C)(\Q)= \Z\cdot [Q_C-O]$, where $Q_C=(-1:1:1)$.\\
Group structure of $G\subset J_0(61)(\Q)$: $G= \Z\cdot D_1\oplus \Z/5\Z\cdot D_{\mathrm{tor}}$, where $D_{\mathrm{tor}}=[(1:0:1:1)-(1:0:0:0)]$ and $D_1=[P+\overline{P}-(1:1:1:1)-(1:0:1:1)]=\rho^*[Q_C-O]$ for 
\[
P:=\left(0: \frac12(\sqrt{-3} - 1): 1: 1\right)\in X_0(61)(\Q(\sqrt{-3}))
\]
satisfying $\rho(P)=Q_C$.\\
Primes used in sieve: 7.\\
There are {\bf no} quadratic points on $X_0(61)$ not coming from $X_0^+(61)(\Q)$ via $\rho: X_0(61)\to X_0^+(61)$ and {\bf no} non-cuspidal rational points. In particular, all quadratic points correspond to $\Q$-curves.

\subsection{$X_0(57)$}
Model for $X_0(57)$: 
\begin{align*}
&x_0x_2 - x_1^2 + 2x_1x_3 + 2x_1x_4 - 2x_2^2 - 2x_2x_3 + 3x_2x_4 - x_3^2 - 2x_3x_4 - x_4^2=0,\\
&x_0x_3 - x_1x_2 - 2x_1x_4 + 4x_2x_3 - 6x_2x_4 - x_3^2 + 5x_3x_4 - 5x_4^2=0,\\
&x_0x_4 - x_2^2 + x_2x_3 - 2x_2x_4 - 2x_4^2=0.
\end{align*}
Genus of $X_0(57)$: 5.\\
Cusps: $(1:0:0:0:0),(1:1:0:1:0),(3:3:1:2:1),(3:9/2:-1/2:7/2:1)$.\\
$C=X_0(57)^*=X_0(57)/\langle w_3,w_{19}\rangle$: elliptic curve $y^2 + y = x^3 - x^2 - 2x + 2$.\\
Group structure of $J(C)$: $J(C)(\Q)\simeq \Z$, $[Q_C-O]\mapsto 1$, where $Q_C=(2:-2:1)$.\\
Group structure of $G\subset J_0(57)(\Q)$: $G= \Z\cdot D_1\oplus \Z/6\Z\cdot D_{\mathrm{tor},1}\oplus \Z/30\Z\cdot D_{\mathrm{tor},2}$, where $D_{\mathrm{tor},1}=[(1:1:0:1:0)-(1:0:0:0:0)]$, $D_{\mathrm{tor},2}=[(3:3:1:2:1)-(1:0:0:0:0)]$ and $D_1=[\mathrm{Trace}(P)-\sum_{c\in \mathrm{cusps}}c]$ for
\begin{align*}
P:=&(\frac{1}{13}(-9\al^3 - 2\al^2 - 3\al + 34):
    \frac{1}{13}(15\al^3 + 12\al^2 - 21\al + 30):
    \frac{1}{13}(-5\al^3 - 4\al^2 + 7\al + 3):\\
    &\frac{1}{13}(8\al^3 + 9\al^2 - 19\al + 29):1),\;\; \al^2=\frac{1+\sqrt{-3}}{2}
\end{align*}
satisfying $\rho(P)=Q_C$.\\
The following are {\bf all} quadratic points on $X_0(57)$ up to Galois conjugacy. There are {\bf no} non-cuspidal rational points. \\
\begin{table}[h!]
    \centering
    \begin{tabular}{c c c c c c}
        Name & $\theta^2$ & Coordinates & $j$-invariant & CM & $\Q$-curve  \\[0.5ex]  
        \hline \hline  
        $P_1$ & -23 & \tiny{$\left(\frac{1}{32}(-11  \theta + 47): \frac{1}{16}(-7  \theta + 43): \frac18(  \theta - 5): \frac14(-  \theta + 9): 1\right)$} & $\frac{343\theta + 4021}{4}$ & NO & NO \\   
    $P_2$ & -23 & $\left(\frac12(  \theta - 3): \frac12(  \theta + 3): 1: 1: 0\right)$ &$\frac{402878\theta + 2212325}{32}$&NO & NO \\  
    $P_3$ & -23 & $\left(\frac12(-  \theta + 3): -  \theta + 2: -1: \frac12(-  \theta + 1): 1\right)$ & \tiny{$\substack{(152503825515346075337\theta \\ - 1518681643605456439979)\\/288230376151711744}$}& NO & NO \\    
    $P_4$ & -23 & $ \left(\frac18(  \theta - 1): \frac18(  \theta + 11): \frac14(-  \theta - 1): \frac18(  \theta + 15): 1 \right)$ & \tiny{$\substack{ (-56278625425021601\theta\\ - 102516814328210867)\\/1048576}$} & NO & NO\\
    $P_5$ & -3 &  $\left(2:\frac12(-  \theta + 7): 0: \frac12(-  \theta + 5): 1\right)$&-12288000& -27&$w_{19}$ \\
    $P_6$ & -3 &$\left(-  \theta: -  \theta + 3: \frac12(  \theta - 1): -  \theta + 2: 1\right)$ &0&-3 & $w_{19}$ \\
    $P_7$ & -3 &$\left(\frac12(-  \theta + 3): \frac12(-  \theta + 3): \frac12(  \theta - 1): 2: 1\right)$ &54000&-12 &$w_{19},w_{57}$ \\
    $P_8$ & -3 & $\left(\frac12(-5  \theta + 1): \frac12(-5  \theta+ 1): \frac12(  \theta- 3): -  \theta + 1: 1\right)$ &0& -3&$w_{19},w_{57}$  \\
    $P_9$ & -2 &$\left(\frac13(  \theta + 4): \frac13(4  \theta + 7): \frac13(-  \theta - 1): \frac13(2  \theta + 5): 1 \right)$ &8000&-8 & $w_{57}$ \\
    $P_{10}$ & -2 &$\left(\frac12(-  \theta+ 2): 3:\frac12   \theta: \frac12(-  \theta + 6): 1\right)$  &8000& -8&$w_{57}$  \\
    $P_{11}$ & 13 & $\left(\frac13(2  \theta + 14): \frac13(2  \theta + 14): \frac16(-  \theta - 1): \frac16(5  \theta + 29): 1\right)$&\tiny{$\substack{(3387888351672962316333\theta \\ - 12215205167504087643323)/2}$}&NO &$w_{3}$  \\
    $P_{12}$ & 13 &$\left(\frac12(5  \theta + 23): \frac12(-3  \theta - 3): \frac12(  \theta + 3): \frac12(-  \theta + 3): 1\right)$ & $\frac{-17787\theta - 61763}{2}$&NO&$w_3$ 
    \end{tabular}
    \label{table57}
\end{table}
\noindent Primes used in sieve: 11,13.
\subsection{$X_0(65)$}
Model for $X_0(65)$:
\begin{align*}
&x_0 x_2 - x_1^2 + x_1 x_4 - 2 x_2^2 - x_2 x_3 + 3 x_2 x_4 - x_3^2 + 2 x_3 x_4 - 2 x_4^2=0,\\
&x_0 x_3 - x_1 x_2 - 2 x_2^2 - x_2 x_3 + 4 x_2 x_4 - x_3^2 + 2 x_3 x_4 - 2 x_4^2=0,\\
&x_0 x_4 - x_1 x_3 - 2 x_2^2 - 3 x_2 x_3 + 5 x_2 x_4 + 3 x_3 x_4 - 3 x_4^2=0.
\end{align*}
Genus of $X_0(65)$: 5.\\
Cusps: $(1:0:0:0:0),(1:1:1:1:1),(1/3:2/3:2/3:2/3:1),(1/2:1/2:1/2:1/2:1)$.\\
$C=X_0(65)^+$: elliptic curve $y^2 + xy = x^3 - x$ of conductor 65.\\
Group structure of $J(C)$: $J(C)(\Q)= \Z\cdot [Q_C-O]\oplus \Z/2\Z\cdot [(0:0:1)-O]$, where $Q_C=(1:0:1)$.\\
Group structure of $G\subset J_0(65)(\Q)$: $G= \Z\cdot D_1 \oplus \Z/2\Z\cdot (-9D_{\mathrm{tor},1}+2D_{\mathrm{tor},2})\oplus \Z/84\Z\cdot (17D_{\mathrm{tor},1}+13D_{\mathrm{tor},2})$, where $D_{\mathrm{tor},1}=[(1:1:1:1:1)-(1:0:0:0:0)]$, $D_{\mathrm{tor},2}=[(1/3:2/3:2/3:2/3:1)-(1:0:0:0:0)]$ and $D_1=[P+\overline{P}-(1:0:0:0:)-(1:1:1:1:1)]=\rho^*([Q_C-O])$ for
\[
P=\left(0:1:\frac12(1+i):1:1\right)\in X_0(65)(\Q(i))
\]
satisfying $\rho(P)=Q_C$. \\
There are {\bf no} quadratic points on $X_0(65)$ that do not come from $X_0(65)^+(\Q)$ via $\rho: X_0(65)\to X_0(65)^+$ and {\bf no} non-cuspidal rational points. In particular, all quadratic points correspond to $\Q$-curves.\\
Primes used in sieve: 17, 23.

\subsection{$X_0(67)$}
Model for $X_0(67)$:
\begin{align*}
   &x_0 x_2 - x_1^2 + 2 x_1 x_3 + 2 x_1 x_4 - 2 x_2^2 - 2 x_2 x_3 + 3 x_2 x_4 - x_3^2 - 2 x_3 x_4 - x_4^2=0,\\
&x_0 x_3 - x_1 x_2 - 2 x_1 x_4 + 4 x_2 x_3 - 6 x_2 x_4 - x_3^2 + 5 x_3 x_4 - 5 x_4^2=0,\\
&x_0 x_4 - x_2^2 + x_2 x_3 - 2 x_2 x_4 - 2 x_4^2=0.
\end{align*}
Genus of $X_0(67)$: 5.\\
Cusps: $(1:0:0:0:0),(1/2 : 1 : 1/2 : 1/2 : 1)$.\\
$C=X_0(67)^+$: genus 2 hyperelliptic curve $y^2 = x^6 - 2x^5 + x^4 + 2x^3 + 2x^2 + 4x + 1$.\\
Group Structure of $C$: $J(C)(\Q)= \Z\cdot [Q_1-(1:-1:0)]\oplus \Z[Q_2-(1:-1:0)]$, where $Q_1=(1:1:0)$ and $Q_2=(0:1:1)$. \\
Group Structure of $G\subset J_0(67)(\Q)$: $G = \Z\cdot D_1\oplus \Z D_2\oplus \Z/11\Z\cdot D_{\mathrm{tor}}$, where $D_{\mathrm{tor}}:=[(1/2 : 1 : 1/2 : 1/2 :1)-(1:0:0:0:0)]$, $D_1:=[P_7+\overline{P}_7-(1:0:0:0:0)-(1/2 : 1 : 1/2 : 1/2 :1)]$ and $D_2:=[P_1+\overline{P}_1-(1:0:0:0:0)-(1/2 : 1 : 1/2 : 1/2 :1)]$ for $P_1,P_7$ defined in the table below satisfying $\rho(P_1)=Q_1$ and $\rho(P_7)=Q_2$.\\
There are {\bf no} quadratic points on $X_0(67)^+$ that do not come from $X_0(67)^+(\Q)$ via $\rho: X_0(67)\to X_0(67)^+$ and there is one non-cuspidal rational point. \\
Moreover, we deduce that the following table gives a complete list of {\bf all} quadratic points and non-cuspidal rational points on $X_0(67)$ up to Galois conjugacy. 

\begin{table}[h!]
    \centering
    \begin{tabular}{c c c c c c}
        Name & $\theta^2$ & Coordinates & $j$-invariant & CM & $\Q$-curve  \\ [0.5ex] 
        \hline \hline  
        $P_0$ & - & $\left(\frac34 : \frac{7}{12} : \frac{7}{12} : \frac13 : 1\right)$ & \small{-147197952000} & -67 & - \\
        $P_1$ & -2 &  \tiny{$\left( \frac{1}{18}(- \theta + 4): \frac{1}{18}( \theta + 14): \frac{1}{18} (- \theta + 4): \frac19 (- \theta + 4): 1 \right)$} & 8000 & -8 & $w_{67}$ \\    
    $P_2$ & -3 &  \tiny{$\left( 0: \frac16 ( \theta + 3): \frac16 ( \theta + 3): \frac16 ( \theta + 3): 1 \right)$} &54000&-12 & $w_{67}$ \\  
    $P_3$ & -3 &  \tiny{$\left(  \frac{1}{26} (3  \theta + 5): 1:  \frac{1}{26} (- \theta + 7):  \frac{1}{13} ( \theta + 6): 1 \right)$} & -12288000& -27 & $w_{67}$ \\   
    $P_4$ & -3 & \tiny{$\left(  \frac{1}{91} (18  \theta + 22):  \frac{1}{182} (-27  \theta + 149):  \frac{1}{182} (15  \theta + 79):  \frac{1}{182} (-3  \theta + 57): 1\right)$} & 0 & -3 & $w_{67}$\\
    $P_5$ & -7 &   \tiny{$\left(  \frac{1}{16} ( \theta + 5):  \frac{1}{16} (- \theta + 11):  \frac{1}{16} ( \theta + 5):  \frac{1}{16} ( \theta + 5): 1 \right)$}&-3375& -7&$w_{67}$ \\
    $P_6$ & -7 & \tiny{$\left(  \frac{1}{20} (- \theta + 1):  \frac{1}{20} (- \theta + 21):  \frac{1}{20} ( \theta + 7):  \frac{1}{20} (- \theta+ 9): 1 \right)$} &16581375&-28 & $w_{67}$ \\
    $P_7$ & -11 & \tiny{$\left( 0:  \frac{1}{11} ( \theta + 11):  \frac{1}{22} (- \theta + 11):  \frac{1}{22} ( \theta + 11): 1 \right)$} &-32768&-11 &$w_{67}$ \\
    $P_8$ & -43 &  \tiny{$\left(  \frac{1}{106} (- \theta - 13):  \frac{1}{53} (-2  \theta + 27):  \frac{1}{106} (-5  \theta + 41):  \frac{1}{53} (-2  \theta + 27): 1 \right)$} &-884736000& -43&$w_{67}$ 
    \end{tabular}
    \label{table67}
\end{table}
\noindent Primes used in sieve: 73, 59, 53, 31, 19, 5.

\subsection{$X_0(73)$}
Model for $X_0(73)$:
\begin{align*}
&x_0 x_2 - 2 x_1^2 + 2 x_1 x_2 - 2 x_1 x_4 - x_2^2 + 3 x_2 x_3 + 3 x_3^2 - x_4^2=0,\\
&x_0 x_3 - 1/2 x_1 x_2 - x_1 x_3 + 1/2 x_2^2 - 1/2 x_2 x_3 + x_2 x_4 - 4 x_3^2 + 9/2 x_3 x_4 - 1/2 x_4^2=0,\\
&x_0 x_4 - x_1 x_3 + x_1 x_4 - x_2 x_3 - 5 x_3^2 + 4 x_3 x_4=0.
\end{align*}
Genus of $X_0(73)$: 5.\\
Cusps: $(1:0:0:0:0),(1:1:1:0:0)$.\\
$C=X_0(73)^+$: genus 2 hyperelliptic curve $y^2 = x^6 + 2x^5 + x^4 + 6x^3 + 2x^2 - 4x + 1$.\\
Group Structure of $J(C)$: $J(C)(\Q) = \Z\cdot [Q_1-(1:1:0)]\oplus \Z\cdot [Q_2-(1:1:0)]$, where $Q_1:=(0:-1:1)$ and $Q_2:=(0:1:1)$. \\
Group Structure of $G\subset J_0(73)(\Q)$: $G= \Z\cdot D_1\oplus \Z\cdot D_2 \oplus \Z/6\Z\cdot D_{\mathrm{tor}}$, where $D_{\mathrm{tor}}=[(1:0:0:0:0)-(1:1:1:0:0)]$, $D_1=[P_3+\overline{P}_3-(1:0:0:0:0)-(1:1:1:0:0)]$, $D_2=[P_6+\overline{P}_6-(1:0:0:0:0)-(1:1:1:0:0)]$ for $P_3$, $P_6$ defined in the table below and satisfying $\rho(P_3)=Q_1$ and $\rho(P_6)=Q_2$. \\
The only quadratic points on $X_0(73)$ that do not come from $X_0(73)^+(\Q)$ via $\rho: X_0(73)\to X_0(73)^+$ are $P_1,P_2$ as defined in the table below. There are {\bf no} non-cuspidal rational points. \\
Moreover, we deduce that the following table gives a complete list of {\bf all} quadratic points on $X_0(73)$ up to Galois conjugacy.

\begin{table}[h!]
    \centering
    \begin{tabular}{c c c c c c}
        Name & $\theta^2$ & Coordinates & $j$-invariant & CM & $\Q$-curve  \\ [0.5ex] 
        \hline \hline  
        $P_1$ & -31 &  \tiny{$\left( \frac{1}{32} (-\theta - 33): \frac{1}{16} (-\theta - 9): \frac{1}{32} (-3 \theta - 35): \frac{1}{32} (\theta + 17): 1 \right)$} & \tiny{$\substack{(-218623729131479023842537441\theta\\- 
    75276530483988147885303471)\\/18889465931478580854784}$} & NO & NO \\  
    $P_2$ & -31 &  \tiny{$\left( \frac{1}{32} (-\theta - 31): \frac{1}{16} (\theta - 17): -\frac{3}{2}: \frac{1}{2}: 1 \right)$} &$\frac{-6561\theta + 1809}{4}$&NO & NO \\ 
    $P_3$ & -19 &  \tiny{$\left( 1/7 (\theta - 10): \frac{1}{14} (\theta - 17): \frac{1}{14} (\theta - 17): \frac{1}{14} (\theta + 11): 1 \right)$ } & -884736& -19 & $w_{73}$ \\   
    $P_4$ & -1 & \tiny{$\left( \frac{1}{5} (\theta - 7): \frac{1}{5} (2 \theta - 4): \frac{1}{5} (3 \theta - 6): \frac{1}{5} (\theta + 3): 1 \right)$} &  287496& -16 & $w_{73}$\\
    $P_5$ & -1 &   \tiny{$\left( \frac{1}{13} (2 \theta - 16): \frac{1}{13} (-8 \theta - 14): \frac{1}{13} (-2 \theta - 23): \frac{1}{13} (2 \theta + 10): 1 \right)$}&1728& -4&$w_{73}$ \\
    $P_6$ & -2 &\tiny{$\left( \frac{1}{6} (-\theta - 8): -1: \frac{1}{6} (-\theta - 8): \frac{1}{6} (-\theta + 4): 1 \right)$} & 8000&-8 & $w_{73}$\\
    $P_7$ & -127 & \tiny{$\left( \frac{1}{32} (-\theta - 47): \frac{1}{176} (-5 \theta - 163): \frac{1}{22} (-\theta - 26): \frac{1}{44} (-\theta + 29): 1 \right)$} &\tiny{$\substack{(14758692270140155157349165\theta\\ + 
    81450017206599109708140525)\\/18889465931478580854784}$}&NO & $w_{73}$ \\
    $P_8$ & -3 & \tiny{$\left( \frac{1}{26} (-9 \theta - 41): \frac{1}{55} (-15 \theta - 51): \frac{1}{26} (-9 \theta - 28): \frac{1}{52} (-9 \theta + 37): 1 \right)$} &0&-3 &$w_{73}$ \\
    $P_9$ & -3 &  \tiny{$\left( -1: \frac{1}{4} (\theta - 3): \frac{1}{4} (\theta - 5): \frac{1}{2}: 1 \right)$} &54000& -12&$w_{73}$ \\
    $P_{10}$ &-1&\tiny{$\left( -1: \frac12(\theta - 3): -2: 1: 1 \right)$} &-12288000&-27&$w_{73}$\\
    $P_{11}$ &-67&\tiny{$\left( \frac{1}{23} (-4 \theta - 43): \frac{1}{46} (-7 \theta - 35): \frac{1}{23} (-3 \theta - 15): \frac{1}{23} (-\theta + 18): 1 \right)$}&-147197952000&-67&$w_{73}$
    \end{tabular}
    \label{table73}
\end{table}
\noindent Primes used in sieve: 43, 67, 41, 17, 37, 13.
\section{The hyperelliptic curve $X_0(37)$}
\label{hypersection}
This curve deserves a special section due to its peculiar nature. A model for $X_0(37)$, computed using the Small Modular Curves package in \texttt{Magma}, is given by
\[
X_0(37):\;y^2 = x^6 + 8x^5 - 20x^4 + 28x^3 - 24x^2 + 
    12x - 4.
\]
It is hyperelliptic of genus 2. We note that Mazur and Swinnerton--Dyer \cite[Section 5]{mazursd} already computed a model for this curve, as well as for $X_0(37)^+$, in 1974 without the aid of computers. They noticed that 37 is the smallest value of $N$ such that $X_0(N)^+$ has positive genus. 

Being hyperelliptic, $X_0(37)$ admits a hyperelliptic involution, mapping $(x,y)\mapsto (x,-y)$. Naturally, one is led to wonder: what is its moduli interpretation? Is it perhaps the Atkin--Lehner involution? There are exactly nineteen values of $N$ such that $X_0(N)$ is hyperelliptic, as was discovered by Ogg \cite{ogg}. Only in three of those nineteen cases, the hyperelliptic involution is not an Atkin--Lehner involution. For two of these, $N=40$ and $N=48$, there are different matrices in $\mathrm{SL}_2(\R)$ inducing the hyperelliptic involution. Not for $N=37$, however.

Lehner and Newman computed in 1964 in a page of corrections to their paper \cite{lehner} that the hyperelliptic involution of $X_0(37)$ is not induced by any automorphism of the complex upper half plane. Instead, the two cusps $(1:1:0)$ and $(1:1:1)$ (in the closure in $\P(1,3,1)$ of the above model) are mapped to the non-cuspidal rational points $(1:-1:0)$ and $(1:-1:1)$ respectively. We now know by \cite{mazur0} that in fact these are the only non-cuspidal rational points. This observation of Lehner and Newman means that any moduli interpretation of the hyperelliptic involution must include an interpretation of the cusps. Such an interpretation in terms of ``generalized elliptic curves" was given by Deligne and Rapoport \cite{deligne}, but we do not attempt to describe the hyperelliptic involution in those terms here.

Define $\infty_+:=(1:1:0)$ and $\infty_-:=(1:-1:0)$ and let $i$ be the hyperelliptic involution on $X_0(37)$. We then have three non-trivial involutions $i,w_{37}$ and $i\circ w_{37}=w_{37}\circ i$. In fact, we compute in \texttt{Magma} that $\mathrm{Aut}_{\Q}(X_0(37))=\langle i,w_{37}\rangle\simeq \Z/2\Z\times \Z/2\Z$. We compute $w_{37}$ to be
\[
w:=w_{37}: X_0(37)\to X_0(37), \; (x:y:z)\mapsto (x:y:x-z).
\]
We refer to the four mentioned rational points as $\infty_+$, $\infty_-$, $w(\infty_+)$ and $w(\infty_-)$. The fact that $i\neq w_{37}$ has dramatic consequences for the abundance of quadratic points on $X_0(37)$. The hyperelliptic covering map $x: X_0(37)\to \P^1$ is one source of infinitely many quadratic points. Unlike in the cases $N=40$ and $N=48$ (see \cite{bruin}), here the quotient $\rho: X_0(37)\to X_0(37)^+:=X_0(37)/\langle w_{37}\rangle$ is an elliptic curve of rank 1. A model is given by
\[
X_0(37)^+: \;y^2 + y = x^3 - x
\]
and its Mordell-Weil group is the free abelian group generated by the point $Q_1:=(0:-1:1)$ (choosing $O=\rho(\infty_+)=(0:1:0)$).  
Hence this bielliptic quotient $X_0(37)^+$ is a second source of infinitely many quadratic points. Moreover, we also have the quotient map $\pi: X_0(37)\to E:=X_0(37)/\langle i\circ w_{37}\rangle$. A model for $E$ is
\[
E: \;y^2 + y = x^3 + x^2 - 23x - 50.
\]
This elliptic curve has Mordell-Weil group $\Z/3\Z$, generated by the point $Q_E:=(8:18:1)$ (choosing $O=\pi(\infty_+)=(0:1:0)$). Another source of points, albeit just three. In this Section we describe how these different sources fit together to make up all points on $X_0(37)^{(2)}(\Q)$.

Consider the map
\[
\kappa: X_0(37)^{(2)}(\Q)\to J_0(37)(\Q), \; \{P,Q\}\mapsto [P+Q - \infty_+-\infty_-].
\]
\begin{lemma}
Let $X$ be any hyperelliptic curve of genus 2 with involution $i$. Then an effective degree 2 divisor $D$ is linearly equivalent to the canonical divisor $K$ if and only if $D$ is of the form $P+i(P)$ for some point $P$ on $X$.
\end{lemma}
\begin{proof}
This is well-known. Riemann--Roch implies that $\ell(K-P-i(P))=\ell(P+i(P))-1$, and $\deg (K-P-i(P))=0$. Post-composed with an automorphism of $\P^1$, the hyperelliptic covering map $X\to \P^1$ defines a non-constant function in the Riemann--Roch space $L(P+i(P))$, so that $K=[P+i(P)]$.
Conversely, suppose that $P+Q\sim R + i(R)$ for points $P,Q,R$ on $X$. Note that $\ell(R+i(R))=2$ by Riemann--Roch, hence $L(R+i(R))=\{\al+\beta f : \al,\beta \in \overline{\Q}\}$, where $f$ is the hyperelliptic covering map with poles in $R,i(R)$. This shows that $R+i(R)$ is only linearly equivalent to degree 2 divisors of the form $P+i(P)$. 
\end{proof}
\begin{lemma}
Let $X$ be a hyperelliptic genus 2 curve with Jacobian $J$, hyperelliptic involution $i$ and $\kappa: X^{(2)}\to J, \;P\mapsto [P-\infty-i(\infty)]$ for some point $\infty$ on $X$. The fibre $\kappa^{-1}(\{0\})$ consists of the effective divisors linearly equivalent to the canonical divisor $K$. Moreover, $\kappa$ is a bijection away from this fibre.
\end{lemma}
\begin{proof}
This is again a well-known consequence of Riemann--Roch.
\end{proof}
%
 %
\begin{proposition}
Let $P_1=(1/2(\sqrt{-3} + 1): -1: 1)\in X_0(37)(\overline{\Q})$.  The Jacobian $J_0(37)(\Q)= \Z\cdot D\oplus\Z/3\Z D_{\mathrm{tor}}$, where $D=[P_1+\overline{P_1}-\infty_+-\infty_-]=\rho^*([Q_1-O])$ and $D_{\mathrm{tor}}=[\infty_++w(\infty_-)-\infty-\infty_-]=[\infty_+-w(\infty_+)]=\pi^*[Q_E-O]$.  
\end{proposition}
\begin{proof}
As $X_0(37)$ is hyperelliptic of genus 2, its Mordell--Weil group can be determined using height bounds and an algorithm of Stoll \cite{stoll}. As 37 is prime, we know by Mazur's \cite{mazur1} proof of Ogg's conjecture that the torsion subgroup equals the rational cuspidal subgroup, and we check that the difference of the cusps has order 3, equals $\pi^*([Q_E-O])$ and is of the form $\kappa(\{\infty_+,i(w(\infty_+)\})$.
\end{proof}

Note that $\rho(P_1)\neq Q_1$ because the basedivisor of $\kappa$ is not $w_{37}$-invariant. This gives us a description of all points in $X_0(37)^{(2)}(\Q)$. First, there is the set $\mathcal{Q}_i$ of quadratic points $Q_x$ with rational $x$-coordinate $x$. Though easy to describe coordinatewise, we lack a moduli interpretation of  elliptic curves corresponding to these points. Then we have the points $\mathcal{P}_{1,0}:=\{P_1,\overline{P_1}\}$ and $\mathcal{P}_{0,1}:=\{\infty_+,w(\infty_-)\}$. For any $(a,b)\in \Z\times \Z/3\Z \setminus \{(0,0)\}$, there is another point $P_{a,b}\in  X_0(37)^{(2)}(\Q)$ defined as the unique effective degree 2 divisor $P$ such that $P-\infty_+ -\infty_- \sim a\mathcal{P}_{1,0}+b\mathcal{P}_{0,1}-(a+b)(\infty_++\infty_-)$ for any lift of $b$ to $\Z$. To compute $\mathcal{P}_{a,b}$, one simply computes the 1-dimensional Riemann--Roch space $L(a\mathcal{P}_{1,0}+b\mathcal{P}_{0,1}-(a+b-1)(\infty_++\infty_-)).$ We conclude the following. 
\begin{proposition}\label{prop37}
The map 
\[
\P^1(\Q)\sqcup (\Z\times \Z/3\Z \setminus \{(0,0)\})\to X^{(2)}(\Q),\; x\mapsto Q_x, \; (a,b)\mapsto P_{a,b}
\]
is a bijection. Moreover, $w_{37}$ interchanges the two points in $\mathcal{P}\in X_0(37)^{(2)}(\Q)$ if and only if $\mathcal{P}=P_{a,0}$ for some $a\in \Z$, and in that case $\mathcal{P}$ corresponds to a $\Q$-curve.  
\end{proposition}
\begin{remark}If $P$ is a quadratic point such that $\rho(P)\in X_0(37)^+(\Q)$, then $[P+\overline{P}-\infty_+ -\infty_-]$ is not equal to $\rho^*[\rho(P)-O]=[P+\overline{P}-\infty_+ -w(\infty_+)]$. As $[w(\infty_+)-\infty_-]=-3[P_1-\infty_+-\infty_-]$, this amounts to a shift: if $\rho(P)=n[Q_1-O]$ then $\{P,\overline{P}\}=P_{n-3,0}$. This justifies the claim that $P_{a,b}$ is fixed by $w_{37}$ if and only if $b=0$, despite the incompatiblity of the basedivisors. 
\end{remark}

The three points $\{P,Q\}$ in $X_0(37)^{(2)}(\Q)$ such that $Q=i\circ w_{37}(P)$ are $\{\infty_+,w(\infty_-)\}$, $\{w(\infty_+),\infty_-\}$ and $\{R,\overline{R}\}$, where $R=(2:2\sqrt{37}:1)$, so $E(\Q)$ only contributes one pair of genuinely quadratic points. Yet, this is an interesting pair. We see that $i$ also interchanges the two points, so  $\{R,\overline{R}\}=Q_2$. Again, this is no problem, because $[R+\overline{R}-\infty_+-\infty_-]\neq \pi^*([\pi(R)-O])$ due to $\infty_-$ not mapping to $O$ under $\pi$. As both $i$ and $i\circ w$ interchange $R$ and $\overline{R}$, we find that $w_{37}$ fixes $R$ and $\overline{R}$. The corresponding elliptic curve thus has a degree 37 endomorphism and must have CM by the maximal order in $\Q(\sqrt{-37})$. It is no coincidence that the point is defined over $\Q(\sqrt{37})$: this follows from Theorem 4.1 in \cite{advancedsilverman} as the Hilbert class field of $\Q(\sqrt{-37})$ is $\Q(\sqrt{-37},\sqrt{37})$. The $j$-invariant of the elliptic curve corresponding to $R$ is $3260047059360000\sqrt{37} + 19830091900536000$.

\bibliographystyle{plain}
\bibliography{Quad_pts_v3}

\begin{thebibliography}{10}

\bibitem{abramovich}
Dan Abramovich and Joe Harris.
\newblock Abelian varieties and curves in {$W_d(C)$}.
\newblock {\em Compositio Math.}, 78(2):227--238, 1991.

\bibitem{bbbmtv}
J.~Balakrishnan, A.~Best, F.~Bianchi, B.~Lawrence, S.~M\"uller,
  N.~Triantafillou, and J.~Vonk.
\newblock Two recent $p$-adic approaches towards the (effective) mordell
  conjecture.
\newblock {\em arXiv:1910.12755}.

\bibitem{dogra}
Jennifer~S. Balakrishnan and Netan Dogra.
\newblock Quadratic {C}habauty and rational points, {I}: {$p$}-adic heights.
\newblock {\em Duke Math. J.}, 167(11):1981--2038, 2018.
\newblock With an appendix by J. Steffen M\"{u}ller.

\bibitem{tuitman}
Jennifer~S. Balakrishnan, Netan Dogra, Steffen M\"uller, Jan Tuitman, and Jan
  Vonk.
\newblock Explicit {C}habauty--{K}im for the split {C}artan modular curve of
  level 13.
\newblock {\em Ann. of Math.}, 189(3):885--944, 2019.

\bibitem{bars}
Francesc Bars.
\newblock Bielliptic modular curves.
\newblock {\em J. Number Theory}, 76(1):154--165, 1999.

\bibitem{bruinstoll}
Nils Bruin and Michael Stoll.
\newblock The {M}ordell-{W}eil sieve: proving non-existence of rational points
  on curves.
\newblock {\em LMS J. Comput. Math.}, 13:272--306, 2010.

\bibitem{bruin}
Peter Bruin and Filip Najman.
\newblock Hyperelliptic modular curves {$X_0(n)$} and isogenies of elliptic
  curves over quadratic fields.
\newblock {\em LMS J. Comput. Math.}, 18(1):578--602, 2015.

\bibitem{bugeaud}
Yann Bugeaud, Maurice Mignotte, Samir Siksek, Michael Stoll, and Szabolcs
  Tengely.
\newblock Integral points on hyperelliptic curves.
\newblock {\em Algebra Number Theory}, 2(8):859--885, 2008.

\bibitem{chabauty}
Claude Chabauty.
\newblock Sur les points rationnels des vari\'{e}t\'{e}s alg\'{e}briques dont
  l'irr\'{e}gularit\'{e} est sup\'{e}rieure \`a la dimension.
\newblock {\em C. R. Acad. Sci. Paris}, 212:1022--1024, 1941.

\bibitem{coleman}
Robert~F. Coleman.
\newblock Effective {C}habauty.
\newblock {\em Duke Math. J.}, 52(3):765--770, 1985.

\bibitem{cremona}
J.~E. Cremona.
\newblock {\em Algorithms for modular elliptic curves}.
\newblock Cambridge University Press, Cambridge, second edition, 1997.

\bibitem{deligne}
P.~Deligne and M.~Rapoport.
\newblock Les sch\'{e}mas de modules de courbes elliptiques.
\newblock pages 143--316. Lecture Notes in Math., Vol. 349, 1973.

\bibitem{dkss}
M.~Derickx, S.~Kamienny, W.~Stein, and M.~Stoll.
\newblock Torsion points on elliptic curves over number fields of small degree.
\newblock {\em arXiv:1707.00364v1.1}.

\bibitem{drinfeld}
V.~G. Drinfeld.
\newblock Two theorems on modular curves.
\newblock {\em Funkcional. Anal. i Prilo\v{z}en.}, 7(2):83--84, 1973.

\bibitem{faltings}
Gerd Faltings.
\newblock Diophantine approximation on abelian varieties.
\newblock {\em Ann. of Math. (2)}, 133(3):549--576, 1991.

\bibitem{freitas}
Nuno Freitas, Bao~V. Le~Hung, and Samir Siksek.
\newblock Elliptic curves over real quadratic fields are modular.
\newblock {\em Invent. Math.}, 201(1):159--206, 2015.

\bibitem{harris}
Joe Harris and Joe Silverman.
\newblock Bielliptic curves and symmetric products.
\newblock {\em Proc. Amer. Math. Soc.}, 112(2):347--356, 1991.

\bibitem{hartshorne}
Robin Hartshorne.
\newblock {\em Algebraic geometry}.
\newblock Springer-Verlag, New York-Heidelberg, 1977.
\newblock Graduate Texts in Mathematics, No. 52.

\bibitem{hindry}
Marc Hindry and Joseph~H. Silverman.
\newblock {\em Diophantine geometry}, volume 201 of {\em Graduate Texts in
  Mathematics}.
\newblock Springer-Verlag, New York, 2000.
\newblock An introduction.

\bibitem{kamienny}
S.~Kamienny.
\newblock Torsion points on elliptic curves and {$q$}-coefficients of modular
  forms.
\newblock {\em Invent. Math.}, 109(2):221--229, 1992.

\bibitem{kamienny2}
S.~Kamienny.
\newblock Torsion points on elliptic curves and {$q$}-coefficients of modular
  forms.
\newblock {\em Invent. Math.}, 109(2):221--229, 1992.

\bibitem{kenku}
M.~A. Kenku.
\newblock On the modular curves {$X_{0}(125)$}, {$X_{1}(25)$} and
  {$X_{1}(49)$}.
\newblock {\em J. London Math. Soc. (2)}, 23(3):415--427, 1981.

\bibitem{kim1}
Minhyong Kim.
\newblock The motivic fundamental group of {$ \P^1\setminus\{0,1,\infty\}$} and
  the theorem of {S}iegel.
\newblock {\em Invent. Math.}, 161(3):629--656, 2005.

\bibitem{kim2}
Minhyong Kim.
\newblock The unipotent {A}lbanese map and {S}elmer varieties for curves.
\newblock {\em Publ. Res. Inst. Math. Sci.}, 45(1):89--133, 2009.

\bibitem{kolyvagin}
V.~A. Kolyvagin and D.~Yu. Logach\"{e}v.
\newblock Finiteness of the {S}hafarevich-{T}ate group and the group of
  rational points for some modular abelian varieties.
\newblock {\em Algebra i Analiz}, 1(5):171--196, 1989.

\bibitem{lehner}
J.~Lehner and M.~Newman.
\newblock Weierstrass points of {$\Gamma _{0}\,(n)$}.
\newblock {\em Ann. of Math. (2)}, 79:360--368, 1964.

\bibitem{oggconj2}
San Ling.
\newblock On the {$\bold Q$}-rational cuspidal subgroup and the component group
  of {$J_0(p^r)$}.
\newblock {\em Israel J. Math.}, 99:29--54, 1997.

\bibitem{manin}
Ju.~I. Manin.
\newblock Parabolic points and zeta functions of modular curves.
\newblock {\em Izv. Akad. Nauk SSSR Ser. Mat.}, 36:19--66, 1972.

\bibitem{mazur1}
B.~Mazur.
\newblock Modular curves and the {E}isenstein ideal.
\newblock {\em Inst. Hautes \'{E}tudes Sci. Publ. Math.}, (47):33--186 (1978),
  1977.

\bibitem{mazur5}
B.~Mazur.
\newblock Modular curves and the {E}isenstein ideal.
\newblock {\em Inst. Hautes \'{E}tudes Sci. Publ. Math.}, (47):33--186 (1978),
  1977.
\newblock With an appendix by Mazur and M. Rapoport.

\bibitem{mazur0}
B.~Mazur.
\newblock Rational isogenies of prime degree (with an appendix by {D}.
  {G}oldfeld).
\newblock {\em Invent. Math.}, 44(2):129--162, 1978.

\bibitem{mazursd}
B.~Mazur and P.~Swinnerton-Dyer.
\newblock Arithmetic of {W}eil curves.
\newblock {\em Invent. Math.}, 25:1--61, 1974.

\bibitem{poonen}
William McCallum and Bjorn Poonen.
\newblock The method of {C}habauty and {C}oleman.
\newblock In {\em Explicit methods in number theory}, volume~36 of {\em Panor.
  Synth\`eses}, pages 99--117. Soc. Math. France, Paris, 2012.

\bibitem{merel}
Lo\"{\i}c Merel.
\newblock Bornes pour la torsion des courbes elliptiques sur les corps de
  nombres.
\newblock {\em Invent. Math.}, 124(1-3):437--449, 1996.

\bibitem{ogg}
Andrew~P. Ogg.
\newblock Hyperelliptic modular curves.
\newblock {\em Bull. Soc. Math. France}, 102:449--462, 1974.

\bibitem{oggconj3}
Masami Ohta.
\newblock Eisenstein ideals and the rational torsion subgroups of modular
  {J}acobian varieties {II}.
\newblock {\em Tokyo J. Math.}, 37(2):273--318, 2014.

\bibitem{ozman}
Ekin Ozman and Samir Siksek.
\newblock Quadratic points on modular curves.
\newblock {\em Math. Comp.}, 88(319):2461--2484, 2019.

\bibitem{parent}
Pierre Parent.
\newblock Bornes effectives pour la torsion des courbes elliptiques sur les
  corps de nombres.
\newblock {\em J. Reine Angew. Math.}, 506:85--116, 1999.

\bibitem{oggconj1}
Yuan Ren.
\newblock Rational torsion subgroups of modular {J}acobian varieties.
\newblock {\em J. Number Theory}, 190:169--186, 2018.

\bibitem{ribet}
Kenneth~A. Ribet.
\newblock Abelian varieties over {$\bf Q$} and modular forms.
\newblock In {\em Modular curves and abelian varieties}, volume 224 of {\em
  Progr. Math.}, pages 241--261. Birkh\"{a}user, Basel, 2004.

\bibitem{siksek}
Samir Siksek.
\newblock Chabauty for symmetric powers of curves.
\newblock {\em Algebra Number Theory}, 3(2):209--236, 2009.

\bibitem{advancedsilverman}
Joseph~H. Silverman.
\newblock {\em Advanced topics in the arithmetic of elliptic curves}, volume
  151 of {\em Graduate Texts in Mathematics}.
\newblock Springer-Verlag, New York, 1994.

\bibitem{stein2}
William Stein.
\newblock {\em Modular forms, a computational approach}, volume~79 of {\em
  Graduate Studies in Mathematics}.
\newblock American Mathematical Society, Providence, RI, 2007.
\newblock With an appendix by Paul E. Gunnells.

\bibitem{stein}
William~Arthur Stein.
\newblock {\em Explicit approaches to modular abelian varieties}.
\newblock ProQuest LLC, Ann Arbor, MI, 2000.
\newblock Thesis (Ph.D.)--University of California, Berkeley.

\bibitem{stoll}
Michael Stoll.
\newblock On the height constant for curves of genus two. {II}.
\newblock {\em Acta Arith.}, 104(2):165--182, 2002.

\bibitem{oggconj4}
Hwajong Yoo.
\newblock The index of an {E}isenstein ideal and multiplicity one.
\newblock {\em Math. Z.}, 282(3-4):1097--1116, 2016.

\end{thebibliography}

 \vspace{.5cm}
  \textsc{Mathematics Institute, University of Warwick, CV4 7AL, United Kingdom}
 
 \emph{E-mail address}: \texttt{j.box@warwick.ac.uk}
 \end{document}